\documentclass[10pt,reqno]{amsart}
\usepackage{amsmath,mathrsfs,amssymb,mathtools,xcolor,graphicx,bm}
\definecolor{lightblue}{rgb}{0.68, 0.85, 0.9}
\usepackage[colorlinks=true, citecolor=red, linkcolor=blue]{hyperref}
\usepackage{pgfplots}
\usepgfplotslibrary{patchplots}

\usepackage{stmaryrd}
\usepackage[backend=biber, style=numeric, sorting=nyt,maxnames=99,date=year]{biblatex} 
\DeclareFieldFormat[article]{journaltitle}{\textup{#1}} 
\DeclareFieldFormat[article]{volume}{\textbf{#1}}
\DeclareFieldFormat[article]{title}{\textit{#1}} 
\DeclareSymbolFont{stmry}{U}{stmry}{m}{n} %

\usepackage{lmodern}
 
\numberwithin{equation}{section}
\newtheorem{theorem}{{\bf Theorem}}[section]
\newtheorem{proposition}[theorem]{{\bf Proposition}}
\newtheorem{corollary}[theorem]{{\bf Corollary}}
\newtheorem{lemma}[theorem]{{\bf Lemma}}

\newtheorem*{claim*}{\bf Claim}
\newtheorem{definition}[theorem]{Definition}
\theoremstyle{definition} 
\newtheorem{remark}[theorem]{{\bf Remark}}

\newcommand{\bxi}{{\bm{\xi}}}
\newcommand{\e}{{\epsilon}}
\newcommand{\p}{{\partial}}
\renewcommand{\o}{{\Omega}}
\renewcommand{\div}{\operatorname{div}}

\newcommand{\dist}{\mathrm{dist}}
\renewcommand{\leq}{\leqslant}
\renewcommand{\geq}{\geqslant}
\newcommand{\R}{\mathbb{R}}
\newcommand{\nn}{\mathbf{n}}
\newcommand{\uu}{{\mathbf{u}}}

\newcommand{\dd}{{\mathrm{d}}}
\newcommand{\vv}{\mathbf{v}}

\newcommand{\ww}{\mathbf{w}}
\newcommand{\hh}{\mathbf{H}}

\newcommand{\E}{E_\e[\ue|\Gamma]}
\newcommand{\B}{B_\e[\ue|\Gamma]}
\newcommand{\pe}{\psi_\e}
\newcommand{\hn}{{\mathcal{H}^{1}}}
\newcommand{\ue}{{\uu_\e}}
\newcommand{\uek}{{\uu_{\e_k}}}
\newcommand{\nne}{{\mathbf{n}_\e}}
\newcommand{\no}{{\mathbf{n}_{\p\o}}}
\newcommand{\dg}{\dd_{\widetilde{\Gamma}}}
\newcommand{\dge}{(\frac{\dd_{\widetilde{\Gamma}}}{\e})}
\newcommand{\rk}{{\mathbb{R}^k}}
\newcommand{\hnn}{\mathcal{H}^{n-1}}
\newcommand{\otn}{{\o_{t}^{-}}}
\newcommand{\otp}{{\o_{t}^{+}}}
\newcommand{\otknp}{{\o_{t}^{k,\pm}}}
\newcommand{\oted}{\o_{t}^{\e}}
\newcommand{\otep}{{\o_t^{\e,+}}}
\newcommand{\oten}{{\o_t^{\e,-}}}
\newcommand{\var}{{\vartheta}}

\newcommand{\had}{{\frac{\dist_N}{2}}}
\newcommand{\pn}{{P_N}}

\newcommand{\g}{{\Gamma}}
\newcommand{\rt}{{\mathbb{R}^2}}
\newcommand{\esssup}{{\mathrm{ess\,sup}}}
\newcommand{\df}{{\dd_F}}
\newcommand{\cf}{{c_F}}
\newcommand{\dn}{{\dd_N}}
\newcommand{\nj}{{\bm{\bnu}_j}}
\newcommand{\dprr}{{\delta^{\prime\prime}}}
\newcommand{\bnu}{{\bm{\nu}}}
\newcommand{\br}{\mathbb{R}}
\newcommand{\dpr}{{\delta^{\prime}}}

\newcommand{\uez}{{\uu_{\e,0}}}
\newcommand{\pek}{\psi_{\e_k}}

\def\b{\big}
\def\rn{\mathbb{R}^n}
\def\lef{\left(}
\def\rig{\right)}

\usepackage{pgfplots}
\pgfplotsset{compat=1.18}
\usepackage{tikz-3dplot}

\begin{document}

\author[]{Xingyu Wang}
\email{xingyuwang@sjtu.edu.cn}
\address{School of Mathematical Sciences, Shanghai Jiao Tong University, 200240 Shanghai, China}
\title[]{The vector-valued Allen-Cahn equation with potentials of high-dimensional double-wells under Robin boundary conditions}
\begin{abstract}
     This work investigates the vector-valued Allen-Cahn equation with  potentials of high-dimensional double-wells under Robin boundary conditions.
 We establish local-in-time convergence of solutions to mean curvature flow with a fixed contact angle  $0<\alpha\leq 90^\circ$, for a broad class of boundary energy densities and well-prepared initial data.  The limiting sharp-interface system is derived, comprising  harmonic heat flows in the bulk and minimal pair conditions at phase boundaries. 
The analysis combines the relative entropy method with gradient flow calibrations and weak convergence techniques. These results extend prior works on the analysis of the vector-valued case without boundary effects (Comm. Pure Appl. Math., 78:1199-1247, 2025) and  the scalar-valued case with boundary contact energy (Calc. Var. Partial Differ. Equ., 61:201, 2022).
\end{abstract}


\maketitle
\noindent

\tableofcontents

\section{Introduction}
\subsection{Background and motivation} The study of phase transitions in vector-valued systems with boundary effects is fundamental in materials science and fluid dynamics. The vector-valued Allen-Cahn equation with high-dimensional double-well potentials subject to nonlinear Robin boundary conditions, given by 
\begin{subequations}\label{int0}
    \begin{align}
    \p_{t}\ue&=\Delta \ue-\frac{1}{\e^2}\nabla F(\ue)&&\text{ in }\o\times(0,T),\label{int0 1}\\
    \ue&=\uu_{\e,0}&&\text{ in }\o\times\{0\},\\
    \p_{\no}\ue&=\frac{1}{\e}\nabla\sigma(\ue)&&\text{ on }\p\o\times(0,T),
\end{align}
\end{subequations}
 arises as the $L^2$-gradient flow (accelerated by $\frac{1}{\e}$) of the Ginzburg-Landau energy with boundary contact energy
 \begin{equation}
   E_\e(\ue):=\int_\o\lef\frac{\e}{2}|\nabla\ue|^2+\frac{1}{\e}F(\ue)\rig dx+\int_{\p\o}\sigma(\ue)d\mathcal{H}^{n-1}. \label{int1}
\end{equation}
 Here, $\o\subset\R^n$ is a bounded smooth domain, $\ue:\o\to\rk$ represents a vector-valued order parameter, $F:\rk\to[0,\infty)$ is a high-dimensional double-well potential  that  vanishes exactly on  two disjoint submanifolds $N^\pm\subset\rk$,  $\no$ denotes the inward-pointing unit normal to $\p\o$,  $\sigma:\rk\to[0,\infty)$ is the boundary energy density, and $\hnn$ denotes the $n-1$-dimensional Hausdorff measure. These settings corresponding to current work will be detailed in Section \ref{setting}.

 In the scalar-valued case where $N^\pm=\{a^\pm\}\subset\R$, Cahn \cite{cahn1977critical}  first introduced the Ginzburg-Landau functional~\eqref{int1} to study the wetting phenomena. Under mass constraint and non-negativity conditions, Modica \cite{modica1987gradient} later established the $\Gamma$-convergence of this functional as $\e\to0$, characterizing the limiting functional. Owen and Sternberg \cite{owen1992gradient} subsequently derived  formal asymptotic expansions for the associated gradient flow, i.e. \eqref{int0}, revealing two fundamental features.
 \begin{itemize}
     \item  At the slow $O(\e)$ time scale, as $\e\downarrow0$, the sharp interface $\g_t$  appears and evolves via mean curvature flow (MCF) in $\o$.
     \item At the boundary $\p\o$, $\g_t$ satisfies a fixed $\alpha\in(0,90^\circ)$ contact angle condition governed by the boundary energy density $\sigma$ and Young's law.
 \end{itemize} 

\begin{center}
\begin{minipage}{0.4\textwidth}
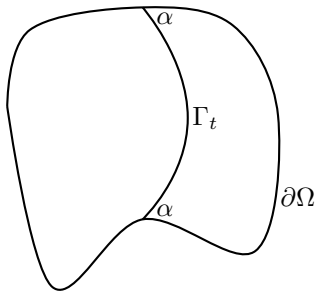
    \textbf{Figure 1.} The interface $\g_t$ meets $\p\o$ at angle $\alpha$.
\end{minipage}
\hfill
\begin{minipage}{0.4\textwidth}
\begin{tikzpicture}[scale=0.6]
\draw[thick]  plot [smooth,tension=0.6] 
coordinates {(-3,1)(-2.5,2.7) (0,3.2) (2,2.8) (3,0.8)
(2.5,-2.2)  (0,-1.5) (-2,-3) (-3,1)};
\draw[thick]  
    plot [smooth, tension=1] 
    coordinates {(0,3.2) (1,0.75) (0,-1.5)};
\node[right] at (0.9,0.75) {$\g_t$};
\node[right] at (0.1,-1.3) {$\alpha$};
\node[right] at (0.1,2.95) {$\alpha$};
\node[right] at (2.85,-1) {$\p\o$};
\end{tikzpicture}
\end{minipage}
\end{center}
 For the case $\alpha = 90^\circ$, corresponding to homogeneous Neumann conditions, Chen  \cite{chen1992generation} established the global-in-time convergence of solutions to sharp-interface limit, utilizing comparison principles and the notion of viscosity solutions.  More recently,  Ables and Moser \cite{abels2019convergence} obtained the local-in-time convergence via the matched asymptotic expansion approach of De Mottoni and Schatzman \cite{de1995geometrical}.  For general contact angles $0<\alpha\leq 90^\circ$, Hensel and Moser \cite{hensel2022convergence} obtained quantitative convergence rates using a relative entropy method and Marshall-Stevens et al. \cite{marshall2024gradient} established  the global-in-time convergence in the framework of varifolds.  Related studies have explored dynamic \cite{abels2022convergence,moser2023convergence}, long-time dynamical behavior \cite{hensel2021bv,kagaya2019convergence,katsoulakis1995generalized,mizuno2015convergence}, and static cases \cite{kagaya2019convergence}.

In the vector-valued setting, equation~\eqref{int0 1} with Neumann boundary conditions was first introduced by Rubinsten et al. \cite{rubinstein1989fast,rubinstein1989reaction} in the context of fast reaction and slow diffusion in chemical reactions.  Their analysis revealed a link between MCF and harmonic heat flows into target manifolds. A significant advance came when Lin et al. \cite{lin2012phase} discovered a \textit{minimal pair condition} (cf.~\eqref{minimal pair} below) for the limiting sharp-interface system by analyzing minimizers of~\eqref{int1} under suitable Dirichlet constraints (without boundary contact energy). The geometric requirement is trivial in the scalar case but non-trivial for a pair of general target manifolds $N^\pm$. Recently, Liu \cite{liu2022phase} extended their results to the dynamical setting under Dirichlet boundary conditions, employing the relative entropy method and weak convergence techniques to obtain local-in-time convergence to MCF and the minimal pair condition. For further related works, see \cite{fei2023matrix} for  $N^\pm=O^\pm(n)$, the $n$-dimensional orthogonal group, and \cite{lin2019harmonic} for the analysis of the limiting sharp-interface system.

Building on these foundations, in this work, we extend the scalar and vector-valued theories.  The aim is twofold: (1) in the vector-valued case, existing rigorous analysis primarily focuses on internal cases (without boundary  contact energy), whereas the case with boundary contact energy is non-trivial and requires intricate analysis due to the high-dimensional nature of vector fields and the coupling between boundaries and interfaces; (2) the special case under Neumann boundary conditions is the original model of the pioneering works by Rubinstein et al. \cite{rubinstein1989fast,rubinstein1989reaction} whose asymptotic analysis in the vector-valued framework remains underexplored, warranting further theoretical investigation.

Our analysis builds upon the relative entropy method (cf. Section \ref{relative rentropy method}) and gradient flow calibrations (cf. Section \ref{section adapted}). 

\begin{itemize}
    \item The relative entropy method was originally developed by Fischer \cite{fischer2020convergence}, with foundations from Jerrard and Smets \cite{jerrard2015motion}, and Fischer and Hensel \cite{fischer2020weak}. This approach has seen significant extensions in recent years, including (1) the matrix-valued generalization by Laux and Liu \cite{laux2021nematic} for studying isotropic-nematic transitions in the Landau-De Gennes model of liquid crystals, (2) Liu's adaptation \cite{liu2024phase} for anisotropic Ginzburg-Landau systems, and (3) the adaptation by Kroemer and Laux \cite{kroemer2025quantitative} for nonlocal Allen-Cahn systems.
    \item The notion of gradient flow calibrations was first introduced by Fischer et al. \cite{fischer2020weak} in their exploration of weak-strong uniqueness for multiphase curvature flows, and subsequently adapted by Hensel and Moser \cite{hensel2022convergence}  to handle boundary contact phenomena. See also \cite{fischer2025local,hensel2022weak,laux2024weak,laux2025diffuse} and references therein for related works.
\end{itemize}

We now articulate the principal challenges of this work. First, the proof of convergence to the mean curvature flow relies on deriving a relative entropy inequality (cf. \eqref{relative entropy inequality 1}) and applying Gr\"onwall's inequality.  While these techniques build on prior works, they require non-trivial adaptations to the setting of this work, refer to Section~\ref{The evolution of the relative entropy}. Second, the derivation of the minimal pair condition presents a key challenge: improving the regularity of the sharp-interface limits \(\uu^\pm\), obtained via weak convergence combined with relative entropy estimates. However, a critical limitation arises because relative entropy estimates only provide
\textit{\(H^1_{\text{loc}}(\o\setminus\Gamma_t)\)-estimates}, leaving the limits ill-defined  on the interface. To overcome this limitation, we leverage the compactness of the space of special functions of bounded variation ($SBV$) to upgrade the regularity of $\uu^\pm$ (cf. Proposition~\ref{regularity}). This uniform bound in $SBV$ is obtained by (1) controlling the discrepancy between the interior reduced boundary of the upper/lower level sets of the quasi-distance function \(\df\) (cf. \eqref{df}) and the interface (cf. Proposition~\ref{level set estimate}); (2) analyzing the relationship between the gradient of \(\uu_\e\) and its tangential projection onto the manifolds \(N^\pm\)(cf. Proposition~\ref{pr core}). Once the regularity is improved, Liu's strategy can be adapted to derive the boundary pair condition.  In particular, to handle boundary terms, we employ tools from geometric measure theory in these analyses, including the precise representative, the reduced boundary, the generalized Gauss-Green formula, and the generalized chain rule, see Section \ref{geometric measure to0ls}. 
This constitutes a significant refinement over \cite{liu2022phase}, where only bulk terms were considered.

\subsection{Settings}\label{setting}
We consider the vector-valued Allen-Cahn equation with a nonlinear Robin boundary condition, given by
\begin{subequations}\label{allen-cahn equation}
    \begin{align}
    \p_{t}\ue&=\Delta \ue-\frac{1}{\e^2}\nabla F(\ue)&&\text{ in }\o\times(0,T),\\
    \ue&=\uu_{\e,0}&&\text{ in }\o\times\{0\},\\
    \p_{\no}\ue&=\frac{1}{\e}\nabla\sigma(\ue)&&\text{ on }\p\o\times(0,T).
\end{align}
\end{subequations}
Here, $\o\subset\mathbb{R}^2$ is a bounded and simply connected domain with a smooth and orientable boundary $\p\o$, $\uu_\e:\o\rightarrow\rk$ with fixed $k\geq2$ is a multi-phase order parameter, and $\no$ denotes the inward-pointing unit normal to $\p\o$. We emphasize that our analysis is  restricted to two-dimensional domains $\o$, see Remark \ref{remark}~(3) for the reason.


We define a pair of disjoint, smooth,  compact, connected submanifolds without boundary,
\begin{equation}
    N = N^+\cup N^- \subset \mathbb{R}^k,
\end{equation}
which represent the two pure phases of the system. Their separation is quantified by the minimal distance
 \begin{equation}
     \dist_N=\min\{|\uu^+-\uu^-|_\rk:\uu^\pm\in N^\pm\}.
 \end{equation}
For $\uu\in\mathbb{R}^k$, we define
\[
\dd_{N^+}(\uu) := \dist(\uu,N^+),\quad\dd_{N^-}(\uu):=\dist(\uu,N^-),
\]
where $\dist(\cdot,N^\pm)$ denotes the Euclidean distance to the set $N^\pm$. The distance $\dd_N$ to $N$ is then defined by
\begin{equation}\label{distance to N}
    \dd_N(\uu):=\min\{\dd_{N^+}(\uu),\dd_{N^-}(\uu)\}.
\end{equation}
Throughout, we fix a small constant $\delta_N\in(0, \frac{\dist_N}{4})$ such that  the nearest-neighbor projection $\pn:B_{2{\delta_N}}(N)\rightarrow N$ is well-defined and smooth. In the $2\delta_N$-tubular neighborhood,
\begin{equation}
    \dd_N(\uu)=|\uu-\pn(\uu)|\quad\forall \uu\in B_{2{\delta_N}}(N).\label{delta_0}
\end{equation}
By the compactness of $N$, there exists a radius $R_N>0$ such that 
\begin{equation}
    B_{\frac{\dist_N}{2}}(N)\subset\{\uu\in\rk:|\uu|< R_N\}=:B_{R_N}.\label{R_0}
\end{equation}

\begin{center}
\begin{minipage}{0.4\textwidth}
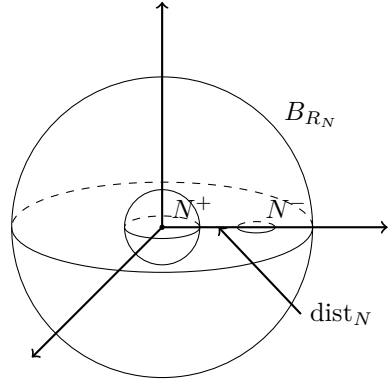
\textbf{Figure 2.} For example, $N^+$ is a two-dimensional sphere and $N^-$ is a one-dimensional circle in $\mathbb{R}^3$.
\end{minipage}
\hfill
\begin{minipage}{0.4\textwidth}
\begin{tikzpicture}
\draw[thick,->] (0,0,0) -- (3,0,0) node[anchor=north east]{};
\draw[thick,->] (0,0,0) -- (0,3,0) node[anchor=north west]{};
\draw[thick,->] (0,0,0) -- (0,0,4.5) node[anchor=south]{};
            \draw[thick,<-] (0.75,0,0) -- (3,0,3);
           \draw[dashed] (0.5,0,0) -- (1,0,0);
           \node[above right] at (0,0,0) {$N^+$};
           \node[above right] at (1.25,0,0){$N^-$};
    \node[right] at (1.5,1.5){$B_{R_N}$};
\draw[dashed] (1.5,0) arc (0:180:0.25 and 0.075);
\draw (1,0) arc (180:360:0.25 and 0.075);  
\node[right] at (3,0,3) {$\dist_N$};

  \draw (0,0) circle (2cm);
  \draw (-2,0) arc (180:360:2 and 0.6);
  \draw[dashed] (2,0) arc (0:180:2 and 0.6);
  \fill[fill=black] (0,0) circle (1pt);
  
   \draw (0,0) circle (0.5cm);
  \draw (-0.5,0) arc (180:360:0.5 and 0.15);
  \draw[dashed] (0.5,0) arc (0:180:0.5 and 0.15);
  \fill[fill=black] (0,0) circle (0.25pt);
\end{tikzpicture}

\end{minipage}
\end{center}

\vspace{2em}

 Next, we specify assumptions on  the potential $F$ and boundary energy density $\sigma$.
 
\vspace{0.5em}
\noindent\textbf{Potential $F$.} 
 The smooth double-well potential $F:\rk\rightarrow[0,\infty)$ vanishes exactly on $N$ and satisfies:
 \begin{itemize}     
     \item[(1)] \textit{Local behavior near $N$}.   $F(\uu)$ depends  only on $\dn(\uu)$ near $N$. Specifically, there exist $R_1>R_N$ and a smooth function $f:[0,\infty)\rightarrow[0,\infty)$ with
     \begin{equation}
        f(s)= \begin{cases}
             s&\text{ if } 0\leq s\leq{\delta_N^2},\\
             2{\delta_N^2}&\text{ if }s\geq2{\delta_N^2},\label{f}
         \end{cases}
     \end{equation}
       such that
     \begin{equation}
         F(\uu)=f\big(\dd_N^2(\uu)\big)\quad\text{ for } |\uu|\leq R_1. \label{F 2}
     \end{equation}
     
     \item [(2)] \textit{Growth at infinity}. There exist a constant $c_1>0$ and $R_2>0$ such that
     \begin{equation}
         F(\uu)\geq c_1|\uu|^2 \quad\text{ for }|\uu|\geq R_2.\label{convexity}
     \end{equation}
   
 \end{itemize}
    The construction \eqref{F 2}  was actually first used to construct weak solutions for harmonic heat flow into a single manifold, see \cite{chen1989existence} and \cite[Section 7.6]{lin2008analysis}. 
  

\vspace{0.5em}
\noindent\textbf{Quasi-distance function.} 
To  describe the boundary energy density $\sigma$, we introduce the quasi-distance function $\dd_F:\rk\rightarrow [0,c_F]$, also known as the phase-field potential, originally developed in \cite{liu2022phase} as a modification of earlier constructions in \cite{fonseca1989gradient,lin2012phase,sternberg1988effect}, defined by 
\begin{equation}
    \dd_F(\uu):=
\left\{
\begin{aligned}
\int_0^{\dd_{N^-}(\uu)}\sqrt{2f(\lambda^2)}d\lambda\quad&\text{ if }\dd_{N^-}(\uu)\leq\tfrac{\dist_N}{2},\\
    \frac{c_F}{2}\quad&\text{ if }\dd_{N}(\uu)>\tfrac {\dist_N}2, \\
    c_F-\int_0^{\dd_{N^+}(\uu)}\sqrt{2f(\lambda^2)}d\lambda\quad&\text{ if }\dd_{N^+}(\uu)\leq\tfrac{\dist_N}{2},
\end{aligned}
\right.\label{df}
\end{equation}
where the surface tension coefficient $c_F$ is given by
\begin{equation}
    c_F:=2\int_0^\had\sqrt{2f(\lambda^2)}d\lambda. \label{cf}
\end{equation}
The function $\dd_F$ is Lipschitz continuous on $\rk$ and satisfies 
\begin{equation}
        \dd_F(\uu)=\begin{cases}
            0 &\text{ iff }\uu\in N^-,\\
            c_F&\text{ iff }\uu\in N^+. 
        \end{cases}
    \end{equation}
The properties of $\dd_F$ will be thoroughly analyzed in Section~\ref{quasi}.

\vspace{0.5em}
\noindent\textbf{Boundary energy density $\sigma$.} The boundary energy density $\sigma\in C^{1,1}\big(\rk;[0,\infty)\big)$  satisfies:
\begin{subequations}
    \begin{align}
       \nabla\sigma(\uu)&=\mathbf{0} \quad&&\text{ for } |\uu|\geq R_N,\label{sigma 1}\\
        \sigma(\uu)&\geq\dd_F(\uu)\cos\alpha \quad&&\text{ for }\uu\in\rk,\label{sigma 2}\\ \sigma(N^-)=\{0\},&\quad \sigma(N^+)=\{c_F\cos\alpha \}, \label{sigma 3}
    \end{align}
\end{subequations}
 where $\alpha\in(0,90^\circ]$ determines the contact angle between  the interface and  $\p\o$.
 
 The decay condition~\eqref{sigma 1} and the growth condition on $F$~\eqref{convexity} ensure the well-posedness of the Allen-Cahn equation~\eqref{allen-cahn equation} (see Lemma~\ref{ue regularity}). The coercivity condition~\eqref{sigma 2} is used to derive the coercivity estimates for the relative entropy functionals defined in~\eqref{E} and~\eqref{B} below. The condition~\eqref{sigma 3} implies Young's law:   $$\sigma(\uu^+)-\sigma(\uu^-)=c_F \cos \alpha \quad \text{ for }\uu^\pm\in N^\pm.$$  These assumptions align with \cite[Eq. (1.9a-b)]{hensel2022convergence}.

\subsection{Main results}

 The first result concerns the local-in-time convergence to the mean curvature flow with a fixed contact angle. 
\begin{theorem}\label{theorem 1}
Let $T\in (0,\infty)$, and let $\bigcup_{t\in [0,T]}\o_t^+{\times}\{t\}$
be a strong solution to mean curvature flow in $\Omega\subset\rt$ 
with a fixed contact angle $\alpha\in (0,90^\circ]$ from Lemma \ref{lemma adapted}. Define 
\begin{equation}
    \Gamma_t:=\overline{\p\o_t^+\cap\o},\quad\Gamma:=\bigcup_{t\in[0,T]}\Gamma_t\times\{t\},\quad\otn:=\o\setminus\overline{\otp}.        
\end{equation}
There exists a sequence of initial data $\uu_{\e,0}\in H^1(\o)$ satisfying
     \begin{align}\label{initial}
      \|\uu_{\e,0}\|_{L^\infty(\o)}\leq R_N, \quad 
       \B(0)+\E(0)\leq C_0\e,   
     \end{align} 
    for some constant $C_0>0$  independent of $\e$, where $\E,\B$ are the relative entropy functionals defined in~\eqref{E}-\eqref{B}. 
    
    For such initial data, for every $\e>0$, there exists a unique weak solution $\uu_\e:\o\times[0,T]\rightarrow\rk$  to the Allen-Cahn equation~\eqref{allen-cahn equation} in the sense of Definition~\ref{weak solution} such that $ \|\uu_{\e}\|_{L^\infty(\o)}\leq R_N$, and 
    \begin{subequations}
        \begin{align}      
            \sup_{t\in[0,T]}\E(t)&\leq C_1\e,\\
            \sup_{t\in[0,T]}\B(t)&\leq C_1\e,\\
             \sup_{t\in[0,T]}\int_\o\b|\dd_F(\ue)-c_F\chi_{\otp}\b|dx&\leq C_1\e^\frac{1}{2},      
        \end{align}
    \end{subequations}
    for some constant   $C_1>0$  independent of $\e$, where $\chi_{\otp}$ is the characteristic function of $\otp$.  
\end{theorem}

The second result concerns the behavior of the limiting system.
\begin{theorem}\label{theorem 2}
   Let $\uu_\e$ be the solution from Theorem~\ref{theorem 1}, and let $\e_k \downarrow 0$ be any sequence. Then, up to a subsequence (not relabeled), there exist limits
    \begin{equation}
        \uu_{\e_k}\xrightarrow{k\to \infty} \uu^\pm\text{ weakly-star in } L^\infty(0,T;H^1_{\mathrm{loc}}(\o_t^\pm)),
    \end{equation}
    where
    \begin{equation}
        \uu^\pm\in L^\infty\big(0,T;H^1(\o_t^\pm;N^\pm
        )\big), \quad \p_t\uu^\pm\in L^2\big(0,T;L^2_{\mathrm{loc}}(\o_t^\pm
        )\big).
    \end{equation}
    The limit $\uu^\pm$ solves the harmonic heat flow into $N^\pm$ in the sense of distributions,
    \begin{equation}\label{heat har}
        \p_t\uu^\pm-\Delta\uu^\pm=A^\pm(\uu^\pm)(\nabla\uu^\pm,\nabla\uu^\pm)\quad\text{ in  } \bigcup_{t\in(0,T)}\o_t^\pm\times\{t\}, 
    \end{equation}
    where $A^\pm$ denotes the second fundamental form of $N^\pm$.
    
    Additionally, for almost every $t\in(0,T)$, the limit $\uu^\pm$ satisfy minimal pair conditions, i.e.
    \begin{equation}
        |\uu^+-\uu^-|_{\rk}(x,t)=\dist_N \quad\text{ for }\hn\text{-a.e. }x\in\Gamma_t\label{minimal pair}.
    \end{equation}
\end{theorem}

\begin{remark}\label{remark}
\leavevmode
\begin{enumerate}
\item These results generalize both the interior case without boundary contact energy studied by Liu~\cite{liu2022phase}, and the scalar-valued case with boundary contact energy analyzed by Hensel and Moser~\cite{hensel2022convergence}.
\item The construction of initial data satisfying condition~\eqref{initial}, is provided in Appendix~\ref{constru}.
\item The analysis is currently restricted to two-dimensional domains $\Omega\subset\rt$. However, if the \emph{boundary adapted gradient flow calibration} (see Lemma~\ref{lemma adapted}) can be extended to higher dimensions, then our methods may be extended.
\item The harmonic heat flow into $N^\pm$ in~\eqref{heat har} holds distributionally: for any test function $\xi \in C^\infty_c\big((0,T) \times \Omega_t^\pm\big)$,
\[
\int \partial_t \uu^\pm \cdot \xi~ dxdt + \int \nabla \uu^\pm : \nabla \xi ~dxdt= \int A^\pm(\uu^\pm)(\nabla \uu^\pm, \nabla \uu^\pm) \cdot \xi ~dxdt.
\]
\end{enumerate}
\end{remark}


\subsection{Notations}

Let $n,k\in\mathbb{N}^+$, we fix the following notations:  
\begin{itemize}
    \item \( \mathbf{u} \cdot \mathbf{v} = \sum_{i=1}^k \uu^i \vv^i \). Dot product of \( \mathbf{u}, \mathbf{v} \in \mathbb{R}^k \).

    \item \( A \mathbf{u} = (A_{ij} \uu^j)_{1 \leq i \leq k} \in \mathbb{R}^n \). Action of  \( A \in \mathbb{R}^{n \times k} \) on  \( \mathbf{u} \in \mathbb{R}^k \).

     \item \( A \cdot \mathbf{u} = A^\top \mathbf{u}\in\mathbb{R}^n \). Dot product between \( A \in \mathbb{R}^{k \times n} \) and \( \mathbf{u} \in \mathbb{R}^k \).
     
    \item \( A:B = \operatorname{tr}(A^\top B) \). Frobenius inner product of \( A, B \in \mathbb{R}^{n \times n} \).
    
    \item \( \mathbf{u} \otimes \mathbf{v} = \mathbf{u} \mathbf{v}^\top \in \mathbb{R}^{n \times n} \). Wedge product of  \( \mathbf{u}, \mathbf{v} \in \mathbb{R}^n \).

    \item \( \div(\mathbf{u}) = \sum_{i=1}^k \frac{\partial u^i}{\partial x_i} \). Divergence of  vector field \( \mathbf{u} = (\uu^1, \dots, \uu^k)^\top \in \mathbb{R}^k \).
    
    \item \( \operatorname{div}(A) = \left( \sum_{j=1}^n \frac{\partial A_{1j}}{\partial x_j}, \dots, \sum_{j=1}^n \frac{\partial A_{nj}}{\partial x_j} \right)^\top \in \mathbb{R}^n \). Divergence of  matrix field \( A = (A_{ij}) \in \mathbb{R}^{n \times n} \).

    \item \( \nabla \mathbf{u} = \left( \frac{\partial \uu^i}{\partial x_j} \right)_{1 \leq i \leq k,1\leq j\leq n} \in \mathbb{R}^{k \times n} \). Gradient of column vector field \( \mathbf{u}  \in \mathbb{R}^n \).

    \item $\p_i=\p_{x_i}$ $(0\leq i\leq n)$, with   $\p_0=\p_t$.

    \item $A\overset{\hnn}{=}B$ implies $\hnn (A\Delta B)=0$, where $A\Delta B=(A\setminus B)\cup(B\setminus A)$ is the symmetric differencce of two sets $A,B$.
    \item  $A\subset_{\hnn}B$ implies $\hnn (B\setminus A)=0$ for two sets $A,B$.
    \item   $\int_{\p\o} u d\hnn$.  Boundary integral via the Sobolev trace operator for $u\in H^1(\o)$. 
    \item    $f:\o\to\mathbb{R}$, $\{f>s\}:=\{x\in\Omega:f(x)>s\}$ for measurable $f:\o\to\mathbb{R}$.
    \item $|A|$.  $n$-Lebesgue measure of  $A\subset\rn$.
    \item  $|A|$.  Frobenius norm of $A\in \mathbb{R}^{n\times k}$, i.e., $|A|:=\sqrt{\operatorname{tr}(A^\top A)}$.
\end{itemize}



\section{Preliminaries}
\subsection{Geometric measure theory}\label{geometric measure to0ls}
In this subsection, we present some results from geometric measure theory, used in Section~\ref{improve the regularity} to establish the regularity of limits \(\uu^\pm\). All results apply to $\o\subset \mathbb{R}^n$ with $n\geq2$.

\vspace{0.5em}
\noindent\textbf{The reduced boundary.} For a set  $E\subset\mathbb{R}^n$
 with finite perimeter, the reduced boundary $\p^*E$ is defined as
the set of $x\in \operatorname{spt} \mu_E$ for which the limit 
\begin{equation}
   \bnu_E:= \lim_{r\downarrow0}\frac{\mu_E\big(B(x,r)\big)}{\big|\mu_E \big|\big(B(x,r)\big)} \quad\text{ exists and belongs to }S^{n-1},\label{measure-theoretic}
\end{equation}
where $\mu_E$ denotes the Gauss-Green measure (cf. \cite[Chapter 15]{maggi2012sets}). 

The following lemma collects  fundamental properties of sets of finite perimeter and the corresponding reduced boundaries.
\begin{lemma}\label{finite c}
    Let $\o\subset\rn$ be a bounded domain with a $C^1$ orientable boundary $\p\o
  $. For a set $E$ of finite perimeter in $\o$, define \(F:=\Omega\setminus E\). Then, $F$ has finite perimeter in \(\Omega\), and the following statements hold:
  
    \begin{subequations}
        \begin{align}
            E \text{ and } F \text{ have finite perimeter in } \br^n\label{perimeter 1},\\             \p^*E\cap\o=\p^*F\cap\o, \quad\bnu_E=-\bnu_F\text{ on this set }  
            \label{perimeter 2},\\
            \no=-\bnu_E\text{ on } \p^*E\cap\p\o,\quad\no=-\bnu_F\text{ on } \p^*F\cap\p\o\label{perimeter 3},\\\label{perimeter 4}
            \p^*E\cap\p^*F\cap\p\o=\emptyset,\\            \hnn\big((\p^*E\cap\p\o)\cup(\p^*F\cap\p\o)\big)=\hnn(\p\o),\label{perimeter 5}
        \end{align}
    \end{subequations}
    where $\no$ denotes the inward-pointing unit
normal to $\p\o$.
\end{lemma}
\begin{proof}
For~\eqref{perimeter 1}, we consider the characteristic function $\chi_{E}:\rn\to\{0,1\}$ as the zero-extension of $\chi_{E}|_{\o}\in BV(\o)$. By \cite[Corollary 3.89]{ambrosio2000functions}, $\chi_E\in BV(\rn)$, which implies $E$ (and thus $F$) has finite perimeter in $\rn$.

Since $E$ and $F$ are sets of finite perimeter in $\rn$, their reduced boundaries $\p^*E$ and $\p^*F$ are well-defined. 

For~\eqref{perimeter 2}, since \(\mu_E = -\mu_F\) as Radon measures in \(\o\) (by the definition of the Gauss-Green measure), the reduced boundaries and normals satisfy \(\p^*E\cap\o = \p^*F\cap\o\) and \(\bnu_E = -\bnu_F\) on this set.  

For~\eqref{perimeter 3}, it follows from the consistency of measure-theoretic outer normals for subsets (see \cite[Exercise 16.6]{maggi2012sets}).

For~\eqref{perimeter 4}, assume $x\in\p^*E\cap\p^*F\cap\p\o$. By comparing the density properties of the topological with those of the reduced boundary (cf. \cite[Corollary 15.8]{maggi2012sets}), we get a contradiction, 
    \begin{equation}
        \frac{1}{2}=\lim_{r\downarrow0}\frac{|(E\cup F)\cap B_r(x)|}{|B_r(x)|}=\lim_{r\downarrow0}\frac{|E\cap B_r(x)|}{|B_r(x)|}+\lim_{r\downarrow0}\frac{|F\cap B_r(x)|}{|B_r(x)|}=1,
    \end{equation}
thus proving~\eqref{perimeter 4}.

   For~\eqref{perimeter 5}, let $\xi\in C^1(\overline{\o},\br^n)$ with $\xi=-\no$ on $\p\o$. By  the divergence theorem,
    \begin{equation}
        \begin{aligned}
           \int_\o\div\xi=\int_{\p\o}-\xi\cdot\no d\hnn =\hnn(\p\o).    
        \end{aligned}
    \end{equation}
    On the other hand, applying the generalized Gauss-Green formula (cf. Lemma \ref{Generalized Gauss-Green formula}) with~\eqref{perimeter 2} and~\eqref{perimeter 3},
    \begin{equation}
    \begin{aligned}
       \int_\o \div\xi =&\int_E\div\xi+\int_F\div\xi\\
         =&\int_{\p^*E}\xi\cdot\bnu_E d\hnn+ \int_{\p^*E}\xi\cdot\bnu_F d\hnn\\
         =&\int_{\p^*E\cap\o}\xi\cdot\bnu_E d\hnn+\int_{\p^*F\cap\o}\xi\cdot\bnu_F d\hnn\\
         &+\int_{\p^*E\cap\p\o}\xi\cdot\bnu_E d\hnn+\int_{\p^*F\cap\p\o}\xi\cdot\bnu_F d\hnn\\
         =&\hnn(\p^*E\cap\p\o)+\hnn(\p^*F\cap\o),
    \end{aligned} 
    \end{equation}
    which completes the proof.
\end{proof}

Next, we define approximate limits and the jump of a Lebesgue measurable function; we refer the reader to \cite[Section 3.6]{ambrosio2000functions} and \cite[Section 3.1]{maggi2023hierarchyplateauproblemsapproximation}. For a Lebesgue measurable function $u:\o\subset\rn\to\br\cup\{\pm\infty\}$, define the approximate upper/lower limits  at $x\in\o$  by
\begin{equation}\label{preciser resprestive}
    u^-(x)=\inf\big\{t\in\br:x\in\{u>t\}^{(0)}\big\},\quad u^+(x)=\sup\big\{t\in\br:x\in\{u<t\}^{(0)}\big\},
\end{equation}
where the set $\{u>t\}^{(0)}$ denotes points of density zero in $\{u>t\}$. We define the measurable function $[u]=u^+-u^-$ as the approximate jump. Let $\Sigma_u=\{x\in\o:[u](x)>0\}$, and define the precise representative $u^*:\o\setminus\Sigma_u\to\br\cup\{\pm\infty\}$ by 
\begin{equation}
    u^*=u^+=u^-\quad\text{ on } \o\setminus\Sigma_u.
\end{equation}
Since $|\Sigma_u|=0$, $u^+,u^-$ and $u^*$ are all Lebesgue representatives of $u$. For $u\in W^{1,1}(\o)$, we have $\hnn(\Sigma_u)=0$ (see \cite[Section 4.1]{ambrosio2000functions}). By the co-area formula, for a.e. $t\in\br$, $\{u<t\}$ has finite perimeter and $|\{u=t\}|=0$. Thus,
\begin{align}
    \o\cap\p^*\{u<t\}=\o\cap\p^*\{u>t\}\quad\text{ for a.e. }t\in\br.
\end{align}
This implies that $x\in  \o\cap\p^*\{u<t\}\subset\{u<t\}^{(\frac{1}{2})}$  with $x\notin\{u<t\}^{(0)}$, we must have $t\geq u^+(x)$; $x\in\o\cap\p^*\{u>t\}\subset\{u>t\}^{(\frac{1}{2})}$ with $x\notin\{u>t\}^{(0)}$, we have $t\leq u^-(x)$. Due to  $\hnn(\Sigma_u)=0$, we conclude that for a.e. $t\in \br$,
\begin{equation}
     \o\cap\p^*\{u<t\}\subset\{t\geq u^+\}\cap\{t\leq u^-(x)\}\overset{\hnn}=\{u^*=t\}.\label{sigma tu}
\end{equation}
On the other hand, from the co-area formula for Sobolev functions (\cite[Theorem 1.1]{Mal2001TheCF}), we have
\begin{equation}
    \int_\o\hnn\big(\o\cap\p^*\{u<t\}\big)dt=\int_\o\hnn\big(\{u^*=t\}\big)dt.\label{sobolev co-area}
\end{equation}
Combining \eqref{sigma tu} and \eqref{sobolev co-area}, we conclude that
\begin{equation}\label{hn eq}
\begin{aligned}
   \o\cap\p^*\{u<t\}\overset{\hnn}{=}\{u^*=t\}\quad\text{ for a.e. }t\in\br.
\end{aligned}
\end{equation}
Building on these definitions and facts, we establish the following lemma, which further explores the properties of sets defined by level functions of a Sobolev function.
\begin{lemma}\label{precise}
  Let $\o\subset\rn$ be a bounded domain with a $C^1$ orientable boundary $\p\o
  $. For $u\in W^{1,1}(\o)$ and  $\mathcal{L}^2$-a.e. $(s,t)\in \br^2$ with $s<t$, define $E^-:=\{u<s\}$, $E^+:=\{u>t\}$ and $E:=E^-\cup E^+$. Then,
\begin{align}
 E^\pm,E\text{ and }E^c \text{ are  sets of finite perimeter in }\o \text{ and } \rn ,\label{per 1}\\ 
 \hnn\big(\o\cap\p^*E^-\cap\p^*E^+\big)=0,
  \label{per 2}\\
   \o\cap\p^*E\overset{\hnn}=\big(\o\cap\p^*E^-\big)\cup\big(\o\cap\p^*E^+\big),\label{per 3}\\
  \bnu=\bnu^-\quad\text{ for $\hnn$-a.e. }x\in\o\cap\p^*E^-,\label{per 4}\\
  \bnu=\bnu^+ \quad\text{ for $\hnn$-a.e. }x\in \o\cap \p^*E^+,\label{per 5}
  \end{align}
  where $\bnu$ is the measure-theoretic outer normal to $E$ and $\bnu^\pm$ is the measure-theoretic outer normal to $E^\pm$.
\end{lemma}

\begin{proof}
     By the co-area formula, for a.e. $s<t$,
     $\{u< s\}$ and $\{u> t\}$ have finite perimeter in $\o$, and by Lemma~\ref{finite c}, also in $\rn$. By \cite[Theorem 16.3]{maggi2012sets}, $E^\pm$, $E$ and $E^c$ are sets of finite perimeter in $\o$ and $\rn$.
   
     Using~\eqref{hn eq}, one obtains
     \begin{equation}
    \o\cap \p^*E^-\overset{\hnn}{=}\{u^*=s\},\quad\o\cap \p^*E^+\overset{\hnn}{=}\{u^*=t\}.
    \end{equation}
    Since $\{u^*=s\}\cap\{u^*=t\}=\emptyset$ for $s\ne t$, we obtain~\eqref{per 2}.
    
     For~\eqref{per 3}, note that $\{u^*= s\}\subset_{\hnn}\{u>t\}^{(0)}={E^+}^{(0)}$ and $\{u^*= t\}\subset_{\hnn} \{u< s\}^{(0)}= {E^-}^{(0)}$. By \cite[Theorem 16.3]{maggi2012sets} for the union of sets with finite perimeter, we get
\begin{equation}
\begin{aligned}
    &\o\cap\p^*E\\
    \overset{\hnn}
    =&\Big(\o\cap {E^+}^{(0)}\cap\p^*E^-\Big) \cup \Big(\o\cap {E^-}^{(0)}\cap\p^*E^+\Big
    )\\
    \overset{\hnn}
    =&\{u^*=s\}\cup\{u^*=t\}\\
    \overset{\hnn}
      =&\big(\o\cap\p^*E^-\big)\cup\big(\o\cap \p^*E^+\big).
\end{aligned}
\end{equation}
Thus,~\eqref{per 3} holds.

Finally,~\eqref{per 4} and~\eqref{per 5} follow from the consistency of measure-theoretic outer normals for subsets: since \(E^\pm \subset E\), the outer normal of $E$ agrees with \(\bnu^\pm\) on their respective reduced boundary \cite[Exercise 16.6]{maggi2012sets}.
\end{proof}

\vspace{0.5em}
\noindent\textbf{$SBV$ function space.}   
We recall the definition of the special bounded variation space \(SBV(\o;\rk)\). As defined in \cite{ambrosio1995new} (see also \cite[Section 2.1]{liu2022phase}), \(SBV(\Omega;\mathbb{R}^k)\) is a subspace of \(BV(\Omega;\mathbb{R}^k)\), consisting of all functions with bounded variation whose Cantor part of the distributional derivatives vanishes. Precisely,
\begin{equation}
  \uu\in SBV(\Omega;\rk)\iff D\uu=\nabla^a \uu \mathcal{L}^n+(\uu^+-\uu^-)\otimes\bnu_\uu\hnn\llcorner J_\uu,
\end{equation}
where $\nabla^a\uu$ is the absolutely continuous part, $J_\uu$  is the  jump discontinuity set with the measure-theoretic outer normal  $\bnu_\uu$, and $(\uu^+,\uu^-)$ are  approximate one-sided limits on  $J_\uu$.

The following classical result from \cite{ambrosio1995new} provides important convergence properties for sequences in \(SBV(\o;\rk)\). 
\begin{lemma}\label{sbv lemma}
    Let $\psi:[0,\infty)\to[0,\infty]$ be a lower semicontinuous, nondecreasing function such that $\lim_{t\to\infty}\frac{\psi(t)}{t}=\infty$. For an open bounded  $\o\subset\mathbb{R}^n$ and sequence  $\{\uu_k\}\subset SBV(\o;\rk)$ with

    \begin{equation}
        \sup_k\int_\o\psi\big(|\nabla^a\uu_k|\big)dx+\sup_k\int\big|\uu_k^+-\uu_k^-\big|d\hnn<\infty.
    \end{equation}
  If $\uu_k\xrightarrow{k\to\infty}\uu$ weakly-star in $BV(\o)$, then
    \begin{itemize}
        \item[(1)] $\uu\in SBV(\o)$.
        \item[(2)] $\nabla^a\uu_k\xrightharpoonup{k\to\infty} \nabla^a\uu$ weakly in $L^1(\o)$.
        \item[(3)] The jump part $\nabla^j\uu_k\xrightharpoonup{k\to\infty}\nabla^j\uu$ weakly-star as Radon measures.
        \item[(4)] For any convex function $\psi$, 
        \begin{equation}  \label{inequality 11}         \int_\o\psi\big(|\nabla^a\uu_k|\big)dx\leq\liminf_{k\to\infty}\int_\o\psi\big(|\nabla^a\uu|\big)dx.
        \end{equation}
    \end{itemize} 
\end{lemma}

\vspace{0.5em}
\noindent\textbf{The generalized Gauss-Green formula.}  
We conclude this subsection with the generalized Gauss-Green formula, which is derived from \cite[Theorem 3.87]{ambrosio2000functions}. 
\begin{lemma}[Generalized Gauss-Green formula]\label{Generalized Gauss-Green formula}
    Suppose $E$ is a set of finite perimeter in a bounded domain $\o\subset\rn$ with a $C^1$ orientable boundary $\p\o$. Then, the generalized Gauss-Green formula holds: for any $\xi\in C^1(\overline{\o};\rn)$ and $\phi\in H^1(\o;\br)$,
    \begin{equation}        \int_E\xi\cdot\nabla\phi+\int_\o\div\xi\phi=\int_{\p^*E\cap\o}\xi\cdot\bnu_E\phi^*d\hnn-\int_{\p^*E\cap\p\o}\xi\cdot\nn_{\p\o} \phi d\hnn,
    \end{equation}
     where $\bnu_E$ denotes measure-theoretic outer unit normal to $E$, $\no$ is the  inward-pointing unit normal to $\p\o$, and $\phi^*$ is the precise representative of $\phi$. 
\end{lemma}




\subsection{Weak solutions to the Allen-Cahn problem }
We begin by defining the concept of weak solutions to the Allen-Cahn problem, building on the framework established in \cite{hensel2022convergence}.
\begin{definition}[\cite{hensel2022convergence}, Definition 5]\label{weak solution}
    Let $T\in(0,\infty)$, and suppose $\uu_{\e,0}$ has finite energy, i.e., $E(\uu_{\e,0})<\infty$. A measurable function $\uu_\e:\o\times[0,T)\rightarrow\rk$ is called a weak solution to~\eqref{allen-cahn equation} if 
    \begin{itemize}
        \item[(1)] Global regularity. $\uu_\e\in H^1\b(0,T;L^2(\o)\b)\cap L^\infty\b(0,T;H^1(\o)\b)$.
        \item[(2)] Initial condition. $\uu_\e(\cdot,0)=\uu_{\e,0} \text{ a.e. in } \o$.
           
        \item[(3)] Variational formulation.  For all $\xi\in C_c(\overline{\o}\times[0,T);\mathbb{R}^k)$, 
        \begin{equation}
        \begin{aligned}
             \int_0^{T}\int_\o\xi\cdot\p_t\uu_\e dxdt
             =&- \int_0^{T}\int_\o\Big(\nabla\xi:\nabla\uu_\e+\xi\cdot DF(\ue)\Big)dxdt\\
             &+\int_0^{T}\int_{\p\o}\xi\cdot\frac{1}{\e}\nabla\sigma(\uu_\e)d\hn dt.
        \end{aligned}            
        \end{equation}
    \end{itemize}
\end{definition}

The existence of weak solutions can be established using a minimizing movement scheme. The higher-order regularity depends on two key factors: the regularity and decay condition of the function \(\sigma\) given in~\eqref{sigma 1}, and the specific form of the potential function $F$. The proof of existence and the derivation of higher-order regularity closely follow the techniques presented in \cite[Appendix A]{hensel2022convergence}, and an outline of the proof is provided in Appendix~\ref{proof lem 2.6}. We state the main results regarding the existence and properties of weak solutions.
\begin{lemma}\label{ue regularity} 
    In the context of Definition~\ref{weak solution}, if the initial datum $u_{\e,0}\in H^1(\o;\rk)$ satisfies
    \begin{equation}
        \|\uu_{\e,0}\|_{L^\infty(\o)}\leq R_N,
    \end{equation}
    where $R_N$ is given in~\eqref{R_0}. Then, for every $\e>0$, there exists a unique weak solution $\ue$ to~\eqref{allen-cahn equation}, satisfying
    \begin{itemize}
        \item[(1)] Uniform boundedness. The $L^\infty$-bound is preserved by the flow, i.e., 
        \begin{equation}
              \|\uu_{\e}\|_{ L^\infty\b(\o\times(0,T)\b)}\leq R_N.
        \end{equation}
       
        \item[(2)] High regularity and  boundary conditions. We have    
        \begin{equation}
              \quad\qquad \uu_\e\in C\b(0,T;H^1(\o)\b)\cap L^2\b(0,T;H^2(\o)\b), \quad\nabla\p_t\uu_\e\in L^2_{\mathrm{loc}}\b(0,T;L^2(\o)\b).
        \end{equation}
        In particular, the PDE holds pointwise a.e. in $\o\times(0,T)$, and the  Robin boundary condition holds in the trace sense, i.e.,
        \begin{equation}
            \p_{\no}\ue(x,t)=\frac{1}{\e}\nabla\sigma\b(\ue(x,t)\b),\label{boundary condition}
        \end{equation}
         for a.e. $t\in(0,T)$ and $\hn$-a.e. $x\in\p\o$.
         \item[(3)] Interior regularity. We have
        \begin{equation}
            \ue\in C_{\mathrm{loc}}^\infty\b(\o\times(0,T)\b).
        \end{equation}
        \item[(4)] Energy dissipation identity. For every $T^\prime\in(0,T)$, it holds
        \begin{equation}
            E_\e(\uu_\e(\cdot,T^\prime))-E_\e(\uu_{\e,0})=-\int_0^{T^\prime}\int_\o\e\left|\p_t\uu_\e\right|^2dxdt,\label{dissip}
        \end{equation}
       where $E_\e(\ue)$ is defined in \eqref{int1}.
    \end{itemize}
    
\end{lemma}

\subsection{Quasi-distance function}\label{quasi}

We commence by defining the following normal:
\begin{equation}\label{bnu}
    \bnu_N(\uu):=\begin{cases}
        \frac{\uu-P_N\uu}{|\uu-P_N\uu|}&\quad\text{ for }\uu\in B_{2\delta_N}(N)\setminus N,\\
        \mathbf{0}&\quad\text{ for } \uu\in N.
    \end{cases} 
\end{equation}
Recall the definition of $\dd_F$ given in~\eqref{df}. The subsequent lemma, as presented in~\cite{liu2022phase}, summarizes its properties.
\begin{lemma}[\cite{liu2022phase}, Lemma 2.5]\label{lemma df}
    The function $\dd_F:\rk\rightarrow [0,c_F]$ is Lipschitz continuous and belongs to $C^1\b(B_{{\delta_N}}(N)\b)$. Moreover,    
       \begin{equation}
        |\nabla\dd_F(\uu)|\leq\sqrt{2F(\uu)} \quad\text{ a.e. }\uu\in \rk, 
        \end{equation}
        \begin{equation}\label{df pro2}
        \dd_F(\uu)=\begin{cases}
            0 &\text{ iff }\uu\in N^-,\\
            c_F&\text{ iff }\uu\in N^+ ,
        \end{cases}
        \end{equation}  
        \begin{align}
        \frac{\nabla\df(\uu)}{|\nabla\df(\uu)|}&=\mp \bnu_N(\uu)&&\text{ if }\uu\in B_{\delta_N}(N^\pm)\setminus N^\pm,\label{df pro3}\\\
            \nabla\dd_F(\uu)&=\mathbf{0}&&\text{ if } \uu\in N.
        \end{align}        
\end{lemma}
Let $\ue$ be  the weak solution obtained from Lemma~\ref{ue regularity}. Following \cite{laux2018convergence} and \cite{liu2022phase}, we consider the restriction of $\dd_F$ to the affine tangent space
\begin{equation}
    T_{x,t}^\ue:=\ue(x,t)+\operatorname{span}\{\p_t\ue(x,t),\p_1\ue(x,t),\p_2\ue(x,t)\}
\end{equation}
denoted by $\dd_F\big|_{ T_{x,t}^\ue}$. By Lemma~\ref{lemma df}, $\dd_F$ is Lipschitz in $\rk$ and $\df(\ue)\in H^1\b(\o\times(0,T)\b)$ . Applying the generalized chain rule of Ambrosio and Dal Maso \cite{ambrosio1990general}, the restriction $ \dd_F\big|_{ T_{x,t}^\ue}$ is differentiable at $\ue(x,t)$ for a.e.  $(x,t)\in\o\times(0,T)$, and
\begin{equation} \label{general}
    \p_i(\dd_F(\ue))=\lim_{h\downarrow0}\frac{\df(\ue+h\p_i\ue)-\df(\ue)}{h} \quad\text{ for }i=0,1,2.
\end{equation}
Let  $A_d\subset\rk$  denote the set of  classically differentiable points of $\df$, and  let $\Pi(x,t):\rk\to T^\ue_{x,t}-\ue(x,t)$ be the orthogonal projection.  Then, the generalized differential 
$ D\dd_F(\ue):\o\times(0,T)\to \rk$ is  defined as follows: for a.e. $(x,t)$ and all $\vv\in\rk$,
 \begin{equation} \label{def general}
     D\dd_F(\ue)(x,t)\cdot\vv=\begin{cases}
       \lim_{h\downarrow0}\frac{\df(\ue+h\Pi(x,t)\vv)-\df(\ue)}{h}&\text{ if  }\ue(x,t)\notin A_d,\\
       \nabla\df(\ue)\cdot\vv&\text{ if }\ue(x,t)\in A_d.
       
     \end{cases}
 \end{equation}

 We distinguish the notation $D\df(\ue)$  for the generalized differential  and $\nabla \df(\ue)$ for the classical derivative.
 The following generalized differential inequality is of crucial importance in the relative entropy estimates.
 \begin{lemma}
    For the generalized differential defined in~\eqref{def general}, there holds
\begin{equation}
    |D\df(\ue)|(x,t)\leq\sqrt{2F(\ue(x,t))}\quad\text{ for a.e. }(x,t).\label{df a.e.}
\end{equation}
 \end{lemma}
 \begin{proof}
     Let $(x,t)\in\o\times(0,T)$ be a point such that $\dd_F\big|_{ T_{x,t}^\ue}$ is differentiable at $\ue(x,t)$. We distinguish three cases. 

     Case 1: \(\dn(\ue(x,t)) > \had\). In this case, the inequality~\eqref{df a.e.} follows directly from the definitions of \(\dd_F\) and \(D \dd_F(\ue)\).
     
     Case 2: \(\dd_N(\ue(x,t)) = \had\). Let $\vv\in T_{x,t}^\ue-\ue(x,t)$ if $\ue(x,t)\notin A_d$, or $\vv\in\rk$ if $\ue(x,t)\in A_d$. For $h>0$, define $\dd_h:=\min\{\dd_N(\ue+h\vv),\had\}$. By the definition of $\df$ and the Lipschitz coefficient of $\dn$ is $1$, we have 
     \begin{equation*}
     \begin{aligned}
         |\df(\ue+h\vv)-\df(\ue)|
         \leq &\bigg|\frac{\cf}{2}-\int_0^{\dd_h}\sqrt{2f(\lambda^2)}d\lambda\bigg|\\
         \leq &\bigg|\int^{\had}_{\dd_h}\sqrt{2f(\lambda^2)}d\lambda\bigg|\\
         \leq& |\dd_N(\ue+h\vv)-\dn(\ue)|\max_{\big[\dd_h,\had\big]} \sqrt{2f(\lambda^2)}\\
         \leq& h|\vv|\max_{\big[\dd_h,\had\big]} \sqrt{2f(\lambda^2)}.
     \end{aligned}             
     \end{equation*}
      Dividing by $h$ and sending $h\downarrow 0$, the continuity of $f$ and $\dn$ yields
     \begin{equation*}
         \big|D\df\big(\ue(x,t)\big)\big|\leq \sqrt{2f\bigg(\Big(\had\Big)^2\bigg)}=\sqrt{2F\big(\ue(x,t)\big)}.
     \end{equation*}

     Case 3: \(\dd_N(\ue(x,t)) < \had\). This case follows similarly to Case 2 with minor modifications.
     
     This completes the proof.
 \end{proof}

\subsection{Boundary adapted gradient flow calibrations} \label{section adapted}
We recall the result of strong solutions for planar mean curvature flow with a constant contact angle \(0 < \alpha \leq 90^\circ\) and \textit{boundary adapted gradient flow calibrations}, which was originally obtained in \cite{hensel2022convergence}. 
The result provides a framework for using the relative entropy approach.
\begin{lemma}\label{lemma adapted}
Let $\Omega\subset\R^2$ be a bounded domain with $C^3$-boundary, $T>0$ and $\alpha\in(0,90^\circ]$. Let ${\o^+}:=\bigcup_{t\in [0,T]} {\o_t^+}\times\{t\}$ be a strong solution to mean curvature flow with constant contact angle $\alpha$ in the following sense.
\begin{itemize}
	\item[(1)] \textit{Evolving regular partition in $\Omega$}. For all $t\in[0,T]$ the set ${\o_t^+}\subset\Omega$ is open and connected with finite perimeter in $\R^2$ such that $\overline{\partial^*{\o_t^+}}=\partial{\o_t^+}$. The interface $\Gamma_t:=\overline{\partial^\ast{\o_t^+}\cap\Omega}$ is a compact, connected, one-dimensional embedded $C^5$-manifold with boundary such that its interior~$\g_t^\circ$ lies in~$\Omega$ and its boundary~$\partial \g_t$ consists of exactly two distinct points which are located on the boundary of the domain, i.e., $\partial \g_t\subset\partial\Omega$. Moreover, there are diffeomorphisms $\Phi(\cdot,t)\colon\R^2\rightarrow\R^2$ of class $C_t^0C_x^5\cap C_t^1C_x^3$ with
	$\Phi(\cdot,0)=\mathrm{Id}$ and
	\[
	\Phi({\o_0^+},t) = {\o_t^+}, \quad \Phi(\g_0,t)= \g_t,
	\quad\quad \Phi(\partial \g_0,t) = \partial \g_t
	\]
	for all $t\in[0,T]$ such that $\Phi\colon \g_0\times[0,T]\rightarrow \g:=\bigcup_{t\in [0,T]} \g_t\times\{t\}$.
	\item[(2)] \textit{Mean curvature flow.} The interface~$\g$ evolves by~MCF in the classical sense.
	\item[(3)] \textit{Contact angle condition.} Let $\nn_{\g_t}$ denote the inner unit normal of $\g_t$ with respect to ${\o_t^+}$ and let $\nn_{\partial\Omega}$ be the inner unit normal of $\partial\Omega$ with respect to $\Omega$. Let $p_0 \in \partial \g_0$ be a boundary point and let $p_t:=\Phi(p_0,t) \in \partial \g_t$. Then 
	\begin{align}\label{eq_angle}
	\nn_{\g_t}|_{p_t} \cdot \nn_{\partial\Omega}|_{p_t} = \cos\alpha
	\end{align}
	for all $t\in[0,T]$ encodes the contact angle condition.
    \end{itemize}
    
 Then, there exists a triple $(\bxi,\hh,\vartheta)$ of maps and 
constants $c\in (0,1)$ and $C>0$ satisfying:
\begin{itemize}
    \item[(1)] \textit{Regularity conditions.}
\begin{subequations}
\begin{align}
\label{eq_regularityXi}
\bxi &\in C^1\big(\overline{\Omega}{\times}[0,T];\mathbb{R}^n\big)
\cap C\big([0,T];C^2_{\mathrm{b}}(\Omega;\mathbb{R}^n)\big),
\\ \label{eq_regularityB}
\hh &\in C\big([0,T];C^1(\overline{\Omega};\mathbb{R}^n)\cap C^2_{\mathrm{b}}(\Omega;\mathbb{R}^n)\big),
\\ \label{eq_regularityWeight}
\vartheta &\in C^1_{\mathrm{b}}\big(\Omega{\times}[0,T]\big) 
\cap C\big(\overline{\Omega}{\times}[0,T];[-1,1]\big).
\end{align}
\end{subequations}
    \item[(2)] \textit{Consistency and calibration conditions.}
    For each~$t\in [0,T]$ the vector field~$\bxi(\cdot,t)$ models an extension of the 
unit normal of~$\g_t$ and the vector field~$\hh(\cdot,t)$ models 
an extension of a velocity vector field of~$\g_t$ in the precise
sense of the conditions
\begin{subequations}
\begin{align}
\label{eq_consistencyProperty}
\bxi(\cdot,t) &= \nn_{\g_t}
\text{ and } \big(\nabla\bxi(\cdot,t)\big)^\top\nn_{\g_t} = 0
&& \text{along } \g_t, 
\\
\label{eq_quadraticLengthConstraint}
|\bxi|(\cdot,t) &\leq 1 {-} c\min\big\{1,\dist^2\big(\cdot,\g_t\big)\big\}
&& \text{in } \Omega,
\end{align}
as well as
\begin{align}\label{eq_calibration1}
|\partial_t\bxi + (\hh\cdot\nabla)\bxi + (\nabla \hh)^\top\bxi|(\cdot,t)
&\leq C\min\big\{1,\dist\big(\cdot,\g_t\big)\big\}
&& \text{in } \Omega,
\\\label{eq_calibration2}
|\bxi\cdot(\partial_t\bxi+(\hh\cdot\nabla)\bxi)|(\cdot,t)
&\leq C\min\big\{1,\dist^2\big(\cdot,\g_t\big)\big\}
&& \text{in } \Omega,
\\ \label{eq_calibrationEvolByMCF}
|\bxi\cdot \hh + \nabla\cdot\bxi|(\cdot,t) 
&\leq C\min\big\{1,\dist\big(\cdot,\g_t\big)\big\}
&& \text{in } \Omega,
\\
|\bxi\cdot(\bxi\cdot\nabla) \hh|(\cdot,t) 
&\leq C\min\big\{1,\dist\big(\cdot,\g_t\big)\big\}
&& \text{in } \Omega,\label{eq_calibration4}
\end{align}
which are accompanied by the (natural) boundary conditions
\begin{align}
\label{eq_boundaryCondXi}
\bxi(\cdot,t)\cdot\no(\cdot) &= \cos \alpha
&&\text{along } \partial \Omega,
\\\label{eq_boundaryCondVelocity}
\hh(\cdot,t)\cdot\no(\cdot) &= 0
&&\text{along } \partial \Omega.
\end{align}
\end{subequations}
    \item[(3)] \textit{Weight function conditions.}
    For all~$t\in [0,T]$ the weight~$\vartheta(\cdot,t)$ models 
a truncated and sufficiently regular ``signed distance'' of~$\g_t$ 
in the sense that 
\begin{subequations}
\begin{align}
\label{eq_weightNegativeInterior}
\vartheta(\cdot,t) &< 0 
&&\text{in the interior of } {\o_t^+},
\\ \label{eq_weightPositiveExterior}
\vartheta(\cdot,t) &> 0 
&&\text{in the exterior of } {\o_t^+},
\\ \label{eq_weightZeroInterface}
\vartheta(\cdot,t) &= 0
&&\text{along } \g_t,
\end{align}
as well as
\begin{align}
\label{2.16d}
\min\{
\dist\big(\cdot,\g_t\big), 1
\big\}&\geq c|\var|(\cdot,t) &&\text{in }\o,\\
\min\{\dist(\cdot,\partial\Omega),
\dist\big(\cdot,\g_t\big), 1
\big\} &\leq C|\vartheta|(\cdot,t)
&&\text{in } \Omega, \label{eq_weightCoercivity}
\\
\label{eq_weightEvol}
|\partial_t\vartheta + (\hh\cdot\nabla)\vartheta|(\cdot,t)
&\leq C|\vartheta|(\cdot,t)
&&\text{in } \Omega. 
\end{align}
\end{subequations}
\end{itemize}

\end{lemma}


\section{Proof of Theorem~\ref{theorem 1}}\label{relative rentropy method}
\subsection{Definitions of the relative entropy functionals}
We begin by introducing the symbol 
\begin{equation}
    \pe:=\dd_F(\ue).
\end{equation}
 Recall the calibration triple \((\bxi,\hh,\var)\) from Lemma~\ref{lemma adapted}. Building upon the frameworks established in \cite{fischer2020weak,fischer2020convergence,hensel2022convergence,laux2021nematic,liu2022phase}, we define the relative entropy functionals \(\E(t)\) and  \(\B(t)\) as follows,
\begin{equation}
\begin{aligned}
    \E(t):=&\int_\o\left(\frac{\e}{2}\left|\nabla\ue(x,t)\right|^2+\frac{1}{\e}F\b(\ue(x,t)\b)-\bxi\cdot\nabla\pe(x,t)\right)dx\\&+\int_{\p\o}\Big(\sigma\b(\ue(x,t)\b)-\dd_F\b(\ue(x,t)\b)\cos\alpha \Big) d\hn(x),\label{E}
\end{aligned}
\end{equation}
and
\begin{equation}
    \B:=\int_\o\b(c_F\chi_{\otp}(x,t)-\pe(x,t)\b)\vartheta(x,t) dx\label{B}.
\end{equation}
Following \cite{liu2022phase}, we define the approximate normal vector and mean curvature:
\begin{align}
 \nn_\e (x,t)&:=\begin{cases}
 \frac{\nabla \psi_\e }{|\nabla \pe|}&\text{ if } |\nabla \pe| (x,t)\neq 0,\\
\mathbf{0}&\text{ otherwise, }
 \end{cases}
\label{normal diff}\\
\hh_\e (x,t)&:=\begin{cases}
-\left(\e  \Delta \uu_\e  -\frac{1}{\e }D  F(\uu_\e  ) \right)\cdot\frac{\nabla \uu_\e  }{\left|\nabla \uu_\e  \right|} &\text{ if } |\nabla \uu_\e|(x,t)\neq 0,\\
\mathbf{0}&\text { otherwise, }
\end{cases}
 \label{mean curvature app}
\end{align}
where the product $DF(\ue)\cdot\nabla\ue$ is understood columnwise. And, we define the projection of each partial derivative $\p_i \ue$ onto the direction of the generalized differential $D\dd_F(\ue)$ (cf.~\eqref{def general}) as
\begin{equation}\label{Pi}
    \Pi_{\ue}\p_i\ue=\begin{cases}
        \lef\p_i\ue\cdot\frac{D\dd_F(\ue)}{|D\dd_F(\ue)|}\rig\frac{D\dd_F(\ue)}{|D\dd_F(\ue)|}&\text{ if } |D\dd_F(\ue)|\ne 0,\\
        \mathbf{0}&\text{ otherwise. } 
    \end{cases}
\end{equation}

We record the following identities.
\begin{lemma}
    For a.e. $(x,t)\in\o\times(0,T)$ and $1\leq i,j\leq 2$, the following identities hold:  
    \begin{subequations}\label{pro}
    \begin{gather}  
    |\nabla\pe| = |\Pi_{\ue}\nabla\ue||D\dd_F(\ue)|, \label{pro1} \\  
    \big(\partial_i \ue - \Pi_{\ue}\partial_i \ue\big) \cdot \Pi_{\ue}\partial_j \ue = 0, \label{pro2} \\  
    \Pi_\ue \partial_i \ue \cdot \Pi_\ue \partial_j \ue = \nn_\e^i \nn_\e^j |\Pi_\ue \nabla\ue|^2. \label{pro3}  
    \end{gather}
    \end{subequations}
    
where $\nn_\e^j$ denotes the j-th component of $\nne$, and $\Pi_\ue\nabla\ue$ denotes the columnwise projection of $\nabla\ue$.

\end{lemma}
\begin{proof}
   \textit{Proof of~\eqref{pro1}.}
    If $|D\df(\ue)|=0$, the generalized chain rule gives $\p_i\pe=D\df(\ue)\cdot\p_i\ue=0$, so~\eqref{pro1} holds trivially. For $|D\df(\ue)|\ne0$, observe that
    $$|\Pi_{\ue}\p_i\ue|^2|D\dd_F(\ue)|^2=\big(\p_i\ue\cdot D\df(\ue)\big)^2=|\p_i\pe|^2,$$
    and summing over $i=1,2$ and using the identity $|\nabla\pe|^2=\sum_{i=1}^2|\p_i\pe|^2$ yield the result.

    \textit{Proof of~\eqref{pro2}.}    
    When $|D\df(\ue)|=0$, both sides vanish due to $\Pi_\ue\p_j\ue=\mathbf{0}$. Otherwise, 
    \begin{equation}
    \begin{aligned}
         \p_i\ue\cdot\Pi_{\ue}\p_j\ue=\frac{\big(\p_i\ue\cdot D\df(\ue)\big)\big(\p_j\ue\cdot D\df(\ue)\big)}{|D\df(\ue)|^2}
         =\Pi_{\ue}\p_i\ue\cdot\Pi_{\ue}\p_j\ue,
    \end{aligned}    
    \end{equation}    
    thus~\eqref{pro2} holds.

    \textit{Proof of~\eqref{pro3}.}
    If $|D\df(\ue)|=0$ or $\p_i\ue=\mathbf{0}$, \eqref{pro3} holds trivial. Otherwise, using $\nabla\pe=|\nabla\pe|\nn_\e$ and the identity~\eqref{pro1}, we compute
    \begin{equation}
        \begin{aligned}
             \Pi_\ue\p_i\ue\cdot \Pi_\ue\p_j\ue=
             \frac{\p_i\pe\p_j\pe}{|\nabla\pe|^2}|\Pi_{\ue}\nabla\ue|^2
             =\nn_\e^i\nn_\e^j|\Pi_{\ue}\nabla\ue|^2,
        \end{aligned}
    \end{equation}
   completing the proof.
\end{proof}
The following coercivity estimates for the relative entropy functionals are adapted from \cite{laux2021nematic,liu2022phase}.
\begin{lemma}\label{le coer}
    There exists a constant $C=C(\Gamma)>0$ independent of $\e,t$ such that for all $t\in[0,T]$,
    \begin{subequations}
        \begin{align}
        \int_{\p\o}\big(\sigma(\ue)-\psi_\e\cos\alpha \big)d\hn\leq
        \E,\label{energy a}\\
            \int_\o\lef\frac{\e}{2}|\nabla\ue|^2+\frac{1}{\e}F(\ue)-|\nabla\psi_\e
            |\rig dx\leq\E,\label{energy b}\\
            \e\int_\o\left|\nabla\ue-\Pi_{\ue}\nabla\ue\right|^2dx\leq2\E,\label{energy c}\\
            \int_\o\lef\e^\frac{1}{2}|\Pi_{\ue}\nabla\ue|-\e^{-\frac{1}{2}}|D\dd_F(\ue)|\rig^2dx\leq2\E,\label{energy d}\\
              \int_\o\lef\frac{\e}{2}|\nabla\ue|^2+\frac{1}{\e}F(\ue)+|\nabla\psi_\e
            |\rig (1-\bxi\cdot\nne)dx\leq4\E,\label{energy e}\\
              \int_\o\lef\frac{\e}{2}|\nabla\ue|^2+\frac{1}{\e}F(\ue)+|\nabla\psi_\e
            |\rig \min(\dd_\Gamma^2,1)dx\leq C\E,\label{energy f}
        \end{align}
    \end{subequations}
    where $\dd_\g(\cdot,t):=\dist(\cdot,\g_t)$.
\end{lemma}
\begin{proof}
 Using the orthogonal decomposition from the identity~\eqref{pro2}, 
 \[
|\nabla\ue|^2=|\Pi_{\ue}\nabla\ue|^2+|\nabla\ue-\Pi_{\ue}\nabla\ue|^2,
\]
and the identity~\eqref{pro1}, we expand  as
      \begin{align*}
          &\frac{\e}{2}|\nabla\ue|^2+\frac{1}{\e}F(\ue)-\bxi\cdot\nabla\pe\\
          =&\frac{\e}{2}|\nabla\ue|^2+\frac{1}{\e}F(\ue)-|\nabla\pe|+|\nabla\pe|(1-\bxi\cdot\nne)\\
          =&\frac{\e}{2}\big|\nabla\ue-\Pi_\ue\nabla\ue\big|^2+|\nabla\pe|(1-\bxi\cdot\nne)\\&+\frac{1}{2}\lef\e^\frac{1}{2}|\Pi_{\ue}\nabla\ue|-\e^{-\frac{1}{2}}|D\dd_F(\ue)|\rig^2.
      \end{align*}
       Integrating over \(\o\) and noting the boundary term in \(\E\) is non-negative, we obtain coercivity estimates \eqref{energy a}-\eqref{energy d}, and  deduce that
\begin{equation}
    \int_\o|\nabla\pe|(1-\bxi\cdot\nne)dx\leq\E.\label{core1}
\end{equation}  
To obtain estimate~\eqref{energy e}, noting that $0 \leq 1 - \bxi \cdot \nn_\e \leq 2$, we multiply the above estimate~\eqref{core1} by $2$ and add  estimate~\eqref{energy b} weighted by $(1 - \bxi \cdot \nn_\e)$. Finally, estimate~\eqref{energy f} is deduced from the geometric bound $1 - \bxi \cdot \nn_\e \geq c \min(\dd_\Gamma^2, 1)$ (cf.~\eqref{eq_quadraticLengthConstraint}) applied to estimate~\eqref{energy e}. This completes the proof.
\end{proof}

\subsection{The evolution of the relative entropy}
\label{The evolution of the relative entropy}

\begin{proposition}\label{proposition relative}
    There exists a constant $C=C(\Gamma)$ independent of $\e,t$ such that  for  a.e. $t\in(0,T)$,
    \begin{equation}\label{relative entropy inequality 1}
    \begin{aligned}
        \frac{d}{dt}\E&+\frac{1}{2\e}\int_\o\big(\e^2|\p_t\ue|^2-|\hh_\e|^2\big) dx+\frac{1}{2\e}\int_\o\big|\hh_\e
        -\e|\nabla\ue|\hh\big|^2dx\\
        &+\frac{1}{2\e}\int_\o\big|\e\p_t\ue-\div\bxi D\dd_F(\ue)\big|^2dx\leq C\E.
    \end{aligned}   
    \end{equation}
   
\end{proposition}
\begin{proof}
We proceed in six steps, utilizing the regularity  of $\ue$ from Lemma~\ref{ue regularity} and the properties of  the calibration triple \((\bxi,\hh,\var)\) from Lemma~\ref{lemma adapted}.

 \textbf{Step $1$. Boundary term reduction.}\\
Using the contact angle condition $\bxi\cdot\no=\cos\alpha$ (cf.~\eqref{eq_boundaryCondXi}) and integration by parts, we obtain
\begin{equation}
  \begin{aligned}
   &-\int_\o\bxi\cdot\nabla\pe dx-\int_{\p\o}\pe\cos\alpha d\hn\\
   =&\int_\o\div\bxi\pe dx+\int_{\p\o}\pe\bxi\cdot\no d\hn-\int_{\p\o}\pe\cos\alpha d\hn\\ \label{cal 111}
   =&\int_\o\div\bxi\pe dx. 
\end{aligned}  
\end{equation}
Differentiating the right-hand side of~\eqref{cal 111} in time and integrating by parts yield
\begin{equation}\label{cal 222}
\begin{aligned}
    \frac{d}{dt} \int_\o\div\bxi\pe dx =&\int_\o\div\bxi\p_t\pe dx+\int_{\o}\p_t\div\bxi\pe dx\\
    =&\int_\o\div\bxi\p_t\pe dx-\int_\o\p_t\bxi\cdot\nabla\pe dx
\end{aligned} 
\end{equation}
where the boundary term vanishes due to $\p_t\bxi\cdot\no=\p_t(\bxi\cdot\no)=0$ along $\p\o$ (cf.~\eqref{eq_boundaryCondVelocity}).

\textbf{Step $2$. Energy time evolution.}\\
 Combining ~\eqref{cal 111}-\eqref{cal 222} with the energy dissipation identity~\eqref{dissip} for $E_\e(\ue)$,  we derive
\begin{equation}
\begin{aligned}
    &\frac{d}{dt}\E\\=&-\e\int_\o|\p_t\ue|^2dx+\int_\o\div\bxi\p_t\pe dx-\int_\o\p_t\bxi\cdot\nabla\pe dx.   
\end{aligned} 
\end{equation}
 Decomposing $\p_t\bxi$ by the transport equations~\eqref{eq_calibration2}-\eqref{eq_calibration1}, and applying the symbol $\nabla\pe=|\nabla\pe|\nn_\e$, we rewrite
\begin{equation}
    \begin{aligned}
        &\frac{d}{dt}\E\\= &-\e\int_\o|\p_t\ue|^2dx+\int_\o\div\bxi\p_t\pe dx \\
    &\underbrace{-\int_\o\left(\p_t\bxi+(\nabla \hh)^\top\bxi+(\hh\cdot\nabla)\bxi\right )\cdot\left( \nne-\bxi\right)|\nabla\pe| dx}_{=:~\mathrm{I}}\\
    &\underbrace{-\int_\o \left(\p_t\bxi+(\hh\cdot\nabla)\bxi\right)\cdot\bxi|\nabla\pe|dx}_{=:~\mathrm{II}}\\
    &-\int_\o (\nabla \hh)^\top\bxi\cdot\bxi|\nabla\pe| dx+\int_\o\left((\nabla \hh)^\top\bxi+(\hh\cdot\nabla)\bxi\right )\cdot\nabla\pe dx\\
    =&\mathrm{I}+\mathrm{II}-\e\int_\o|\p_t\ue|^2dx+\int_\o\div\bxi\p_t\pe dx\\
    &-\int_\o\nabla\hh:\bxi\otimes(\bxi-\nne)|\nabla\pe|dx+\int_\o(\hh\cdot\nabla)\bxi\cdot\nabla\pe dx.\label{cal 1} 
    \end{aligned}
\end{equation}

\textbf{Step $3$. Vector identity for the convective term.}\\
 Following \cite[ Lemma 12]{hensel2022convergence}, we need the following claim to deal with the convective term $\int_\o(\hh\cdot\nabla)\bxi\cdot\nabla\pe dx$.\\
\textit{Claim $1$.} The following identity holds,
    \begin{equation}
        \begin{aligned}
    &\int_\o(\hh\cdot\nabla)\bxi\cdot\nabla\pe dx\\
    =&\int_\o \div\bxi \hh\cdot \nne|\nabla\pe|dx+\int_\o\nabla\hh:\nne\otimes\bxi|\nabla\pe|dx\\
    &-\int_{\p\o}\div_{\p\o}\hh (\no\cdot\bxi) \pe d\hn-\int_\o\div\hh\bxi\cdot\nabla\pe dx,
\end{aligned}\label{claim 1}
    \end{equation}
where $\div_{\p\o}$ denotes the  surface divergence on $\p\o$.
\begin{proof}[Proof of Claim 1]
Using $\div(\bxi\otimes\hh)=(\hh\cdot\nabla)\bxi+\div\hh\bxi$ and $ \div(\hh\otimes\bxi)=(\bxi\cdot\nabla)\hh+\div\bxi\hh$, one has
\begin{align}
    \int_\o(\hh\cdot\nabla)\bxi\cdot\nabla\pe dx=&\int_\o\div(\bxi\otimes\hh)\cdot\nabla\pe dx-\int_\o\div\hh\bxi\cdot\nabla\pe dx,\label{claim a}
\end{align}
and
 \begin{equation}\label{claim b}
     \int_\o \div(\hh\otimes\xi)\cdot\nabla\pe dx=\int_\o \div\bxi \hh\cdot \nne|\nabla\pe|dx+\int_\o\nabla\hh:\nne\otimes\bxi|\nabla\pe|dx.
 \end{equation}
By summing~\eqref{claim a} and~\eqref{claim b}, the verification of~\eqref{claim 1} reduces to proving the following identity:
\begin{equation}\label{claim 1c}
    \int_\o\div(\bxi\otimes\hh)\cdot\nabla\pe dx-\int_\o \div(\hh\otimes\xi)\cdot\nabla\pe dx=\int_{\p\o}\div_{\p\o}\hh (\no\cdot\bxi) \pe d\hn.
\end{equation}
We turn  the bulk error term in the left-hand side of \eqref{claim 1c} to the boundary error. By applying the divergence theorem twice and utilizing the identity $\div\Big(\div(\bxi\otimes\hh)\Big)=\div\Big(\div(\hh\otimes\bxi)\Big)$, we obtain
\begin{equation} \label{claim 11}
   \begin{aligned}
    &\int_\o\div(\bxi\otimes\hh)\cdot\nabla\pe dx\\
    =&-\int_\o \div\Big(\div(\bxi\otimes\hh)\Big)\pe dx-\int_{\p\o} \div(\bxi\otimes\hh)\cdot\no\pe d\hn \\
    =&-\int_\o \div\Big(\div(\hh\otimes\bxi)\Big)\pe dx-\int_{\p\o} \div(\bxi\otimes\hh)\cdot\no\pe d\hn\\
    =&\int_\o\div(\hh\otimes\bxi)\cdot\nabla\pe dx+\int_{\p\o} \div(\hh\otimes\bxi)\cdot\no\pe d\hn\\
    &-\int_{\p\o} \div(\bxi\otimes\hh)\cdot\no\pe d\hn.    
\end{aligned} 
\end{equation}
Again use the identities $\div(\bxi\otimes\hh)=(\hh\cdot\nabla)\bxi+\div\hh\bxi$ and  $ \div(\hh\otimes\bxi)=(\bxi\cdot\nabla)\hh+\div\bxi\hh$ to  expand the boundary error term,
\begin{equation}
\begin{aligned}
    &\int_{\p\o}\div(\hh\otimes\bxi)\cdot\no\pe d\hn-\int_{\p\o}\div(\bxi\otimes\hh)\cdot\no\pe \hn\\
   = &\int_{\p\o}(\bxi\cdot\nabla)\hh\cdot\no\pe d\hn+\underbrace{\int_{\p\o}\div\bxi\hh\cdot\no\pe d\hn}_{=0\text{ due to }\eqref{eq_boundaryCondVelocity}}\\
  & -\int_{\p\o}(\hh\cdot\nabla)\bxi\cdot\no\pe d\hn
   -\int_{\p\o}\div\hh\bxi\cdot\no\pe d\hn\\
   = &\int_{\p\o}\Big((I-\no\otimes\no)\bxi\cdot\nabla\Big)\hh\cdot\no\pe d\hn\\
   &+\int_{\p\o}\Big((\no\otimes\no)\bxi\cdot\nabla\Big)\hh\cdot\no\pe d\hn\\
   &-\int_{\p\o}(\hh\cdot\nabla)\bxi\cdot\no\pe d\hn-\int_{\p\o}\div\hh\bxi\cdot\no\pe d\hn,\\
\end{aligned}\label{cal 1d}
\end{equation}
where in the last step, we have used $$(\bxi\cdot\nabla)\hh=\Big((I-\no\otimes\no)\bxi\cdot\nabla\Big)\hh+\Big((\no\otimes\no)\bxi\cdot\nabla\Big)\hh$$ to decompose this derivative into the tangential and normal parts along $\p\o$. By the fact $\big((I-\no\otimes\no)\bxi\cdot\nabla\big)(\hh\cdot\no)=0$ on $\p\o$, which is due to $(I-\no\otimes\no)\bxi\in T\p\o$ and $\hh\cdot\no=0$ along $\p\o$, we can rewrite the first term appearing in the last step of~\eqref{cal 1d},
\begin{equation}
\begin{aligned}
    \int_{\p\o}\Big((I-\no\otimes\no)&\bxi\cdot\nabla\Big)\hh\cdot\no\pe d\hn\\
    &=-\int_{\p\o}\Big((I-\no\otimes\no)\bxi\cdot\nabla\Big)\no\cdot\hh\pe d\hn. 
\end{aligned}
    \label{cla 1e}
\end{equation}
Due to $\hh\in T\p\o$ and $\bxi\cdot\no=\cos\alpha$ along $\p\o$, we rewrite the third term, 
\begin{equation}
\begin{aligned}
      -\int_{\p\o}(\hh\cdot\nabla)\bxi\cdot\no\pe d\hn
      =\int_{\p\o}(\hh\cdot\nabla)\no\cdot\bxi\pe d\hn.
      \label{cla 1f}
\end{aligned}
\end{equation}
Exploiting the relations
\begin{equation}
    (\bxi\cdot\nabla)\no\cdot\hh=(\hh\cdot\nabla)\no\cdot\bxi,\quad\Big((\no\otimes\no)\bxi\cdot\nabla\Big)\no\cdot\hh=0
\end{equation}
along $\p\o$, we find adding ~\eqref{cla 1e} and~\eqref{cla 1f} gives zero. To the end, by $ \div\hh-\nabla\hh:\no\otimes\no=\div_{\p\o}\hh$,~\eqref{cal 1d} gives
\begin{equation}
\begin{aligned}
     &\int_{\p\o}\div(\hh\otimes\bxi)\cdot\no\pe d\hn-\int_{\p\o}\div(\bxi\otimes\hh)\cdot\no\pe \hn\\=&
     -\int_{\p\o}\div_{\p\o}\hh (\no\cdot\bxi) \pe d\hn.
\end{aligned}  
\end{equation}
Combining this with~\eqref{claim 11}, it follows~\eqref{claim 1c}, which completes the proof.
\end{proof}
Substituting~\eqref{claim 1} into~\eqref{cal 1}, we reorganize the terms as follows,
\begin{equation}\label{cal 2}
   \begin{aligned}
    &\frac{d}{dt}\E\\
    =&\mathrm{I}+\mathrm{II}-\e\int_\o|\p_t\ue|^2dx+\int_\o\div\bxi\p_t\pe dx\\
    &\underbrace{-\int_\o\nabla\hh:(\bxi-\nne)\otimes(\bxi-\nne)|\nabla\pe|dx}_{=:~\mathrm{III}}+\int_\o\nabla\hh:\nne\otimes \nne|\nabla\pe|dx\\
    &+\int_\o \div\bxi \hh\cdot \nne|\nabla\pe|dx-\int_\o\div\hh\bxi\cdot\nabla\pe dx\\
    &-\int_{\p\o}\pe\cos\alpha\div_{\p\o}\hh d\hn.
\end{aligned} 
\end{equation}

\textbf{Step $4$. The energy-momentum identity.}\\
To estimate non-coercive terms, we need a modified energy-momentum identity (cf. \cite[Lemma A.1]{liu2022phase}) including a boundary correction.\\
\textit{Claim $2$.} The following identity holds:
\begin{equation}
     \begin{aligned}
        &\int_\o \hh\cdot\hh_\e|\nabla\ue|dx+\int_\o\div\hh\left(\frac{\e}{2}|\nabla\ue|^2+\frac{1}{\e}F(\ue)\right)dx\\
&-\sum_{i,j}\int_\o{(\nabla\hh)}_{i,j}\e\p_i\ue\cdot\p_j\ue dx+\int_{\p\o}\sigma(\ue)\div_{\p\o}\hh d\hn=0.\\
    \end{aligned}\label{claim 2}
\end{equation}
\begin{proof}[Proof of Claim 2]
Define the energy-momentum tensor $\mathbf{T}_\e$ with components
\begin{equation*}
    {(\mathbf{T_\e})}_{i,j}=\left(\frac{\e}{2}|\nabla\ue|^2+\frac{1}{\e}F(\ue)\right)\delta_{i,j}-\e\p_i\ue\cdot\p_j\ue.
\end{equation*}
 By the definition of $\hh_\e$ in \eqref{mean curvature app}, the tensor satisfies $\div\mathbf{T}_\e = \hh_\e |\nabla\ue|$. To extract the boundary terms, we test this equation against $\hh$ and integrate by parts,
\begin{align*}
    &\int_\o \hh\cdot\hh_\e|\nabla\ue|dx=\int_\o \hh\cdot\div\mathbf{T_\e}dx \\
    =&-\int_\o\nabla\hh:\mathbf{T_\e}dx-\int_{\p\o}\hh^j{(\mathbf{T_\e})}_{i,j}\nn_{\p\o}^i d\hn\\
    =&-\int_\o\div\hh\left(\frac{\e}{2}|\nabla\ue|^2+\frac{1}{\e}F(\ue)\right)dx+\sum_{i,j}\int_\o{(\nabla\hh)}_{i,j}\e\p_i\ue\cdot\p_j\ue dx\\
    &\underbrace{-\int_{\p\o}\hh^j\left(\frac{\e}{2}|\nabla\ue|^2+\frac{1}{\e}F(\ue)\right)\delta_{i,j}\nn_{\p\o}^id\hn}_{=0\text{ due to }\eqref{eq_boundaryCondVelocity}}+\int_{\p\o}\hh^j\e\p_i\uu_\e^k\p_j\uu_\e^k \nn_{\p\o}^id\hn\\
    =&-\int_\o\div\hh\left(\frac{\e}{2}|\nabla\ue|^2+\frac{1}{\e}F(\ue)\right)dx+\sum_{i,j}\int_\o{(\nabla\hh)}_{i,j}\e\p_i\ue\cdot\p_j\ue dx\\
    &-\int_{\p\o}\sigma(\ue)\div_{\p\o}\hh d\hn, 
\end{align*}
where  the last step uses the equality,  
\begin{align*}
    \int_{\p\o}\hh^j\e\p_i\uu_\e^k\p_j\uu_\e^k \nn_{\p\o}^id\hn
    =&\int_{\p\o}\hh^j\e\p_j\uu_\e^k\frac{1}{\e}\p_k\sigma(\ue)d\hn\\
    =&\int_{\p\o}\hh^j\p_j\b(\sigma(\ue)\b)d\hn\\
    =&-\int_{\p\o}\sigma(\ue)\div_{\p\o}\hh d\hn,
\end{align*}
which follows from the boundary condition   $\p_i\uu^k_\e\nn_{\p\o}^i=\frac{1}{\e}\p_k\sigma(\ue)$ in \eqref{boundary condition} and  tangential integration by parts.
\end{proof}

Adding~\eqref{claim 2} into the right-hand side of~\eqref{cal 2}, we obtain

\begin{equation}
  \begin{aligned}
    &\frac{d}{dt}\E\\
    =&\mathrm{I}+\mathrm{II}+\mathrm{III}\\
    &+\int_\o \Big(-\e|\p_t\ue|^2+\div\bxi\p_t\pe+\hh\cdot\hh_\e|\nabla\ue|+\div\bxi \hh\cdot \nne|\nabla\pe|\Big)dx \\
    &+\underbrace{\int_\o\div\hh\left(\frac{\e}{2}|\nabla\ue|^2+\frac{1}{\e}F(\ue)\right)dx-\int_\o\div\hh\bxi\cdot\nabla\pe dx}_{=:~\mathrm{IV}}\\
    &+\underbrace{\int_\o\nabla\hh:\nne\otimes \nne|\nabla\pe|dx-\sum_{i,j}\int_\o{(\nabla\hh)}_{i,j}\e\p_i\ue\cdot\p_j\ue dx}_{=:~\mathrm{V}}\\
    &+\underbrace{\int_{\p\o}\sigma(\ue)\div_{\p\o}\hh d\hn-\int_{\p\o}\pe\cos\alpha\div_{\p\o}\hh d\hn}_{=:~\mathrm{VI}}.
\end{aligned}  \label{cal 3}
\end{equation}

\textbf{Step $5$.  Completing quadratic forms.}\\
We  evaluate the collection of dynamic and convective terms from the second row of~\eqref{cal 3}, and rearrange them via square completion. Specifically, consider 
\begin{align*}
	\notag-&\e  |\p_t \uu_\e  |^2+   \div \bxi  D  \dd_F   (\uu_\e  )\cdot \p_t \uu_\e  
	+ \div \bxi  \hh   \cdot \nabla \psi_\e+ \hh_\e \cdot \hh |\nabla \uu_\e  |
	\\\notag&= -\frac1{2\e } \Big(  |\e  \p_t \uu_\e  |^2 -2\div \bxi  D \df   (\uu_\e  )\cdot \e  \p_t \uu_\e  
	+(\div \bxi )^2 |D \df   (\uu_\e  )|^2 \Big)
	\\\notag&\quad - \frac1{2\e } |\e  \p_t \uu_\e  |^2 + \frac1{2\e }(\div \bxi )^2 |D \df   (\uu_\e  )|^2
	+ \div \bxi  \hh   \cdot \nabla \psi_\e
	\\\notag&\quad - \frac1{2\e } \Big( |\hh_\e |^2 - 2 \e  |\nabla \uu_\e  | \hh_\e \cdot  \hh + \e ^2 |\nabla \uu_\e  |^2 |\hh|^2\Big)
	+ \frac1{2\e } \Big( |\hh_\e |^2 + \e ^2 |\nabla \uu_\e  |^2 |\hh|^2\Big)
	\\&\notag =  -\frac1{2\e } \Big|\e  \p_t \uu_\e  - \div \bxi D \df   (\uu_\e  ) \Big|^2
	- \frac1{2\e } \Big|\hh_\e  - \e  |\nabla \uu_\e  | \hh \Big|^2
	- \frac1{2\e }  |\e  \p_t \uu_\e  |^2 +\frac1{2\e }  |\hh_\e |^2
	\\&\quad + \frac1{2\e } \Big( (\div \bxi )^2  |D  \dd_F   (\uu_\e  )|^2 + 2\e \div \bxi   \nabla \psi_\e  \cdot \hh + |\e  \Pi_{\uu_\e  }\nabla \uu_\e  |^2 |\hh|^2 \Big)
	\\\notag&\quad +\frac\e {2} \left(|\nabla \uu_\e  |^2- | \Pi_{\uu_\e  }\nabla \uu_\e  |^2\right)  |\hh|^2.
	\end{align*}
Substituting this expansion into~\eqref{cal 3}, we arrive at 
\begin{align*}
    \frac{d}{dt}&\E+\frac 1{2\e }\int_\o \left(\e ^2 \left| \p_t \uu_\e    \right|^2-|\hh_\e |^2\right)dx\\&
    +\frac 1{2\e }\int \Big| \e  \p_t \uu_\e    -\div \bxi D \df   (\uu_\e  )  \Big|^2d x+\frac 1{2\e }\int \Big| \hh_\e -\e |\nabla \uu_\e  | \hh \Big|^2dx\\
  =&\mathrm{I}+\mathrm{II}+\mathrm{III}+\mathrm{IV}+\mathrm{V}+\mathrm{VI}\\
    &+\underbrace{\frac 1{2\e } \int \Big| \div \bxi  |D  \dd_F   (\uu_\e  )|\nn_\e  +\e  |\Pi_{\uu_\e  } \nabla \uu_\e  | \hh\Big|^2dx}_{=:~\mathrm{VII}}+\underbrace{\frac\e {2} \int_\o\left(|\nabla \uu_\e  |^2- | \Pi_{\uu_\e  }\nabla \uu_\e  |^2\right)  |\hh|^2 dx}_{=:~\mathrm{VIII}}.    
\end{align*}

\textbf{Step $6$.  Estimation of remainder terms.}\\
We estimate the terms \( \mathrm{I} \)-\( \mathrm{VIII} \) using previously established coercivity estimates in Lemma~\ref{le coer} and the condition \eqref{eq_calibration1}-\eqref{eq_calibration4}.
Terms $\mathrm{I},\mathrm{II}$ and $\mathrm{III}$  are bounded by~\eqref{energy e} combined with $\dd_\g^2\leq C(1-\bxi\cdot\nne)$ and $|\nne-\bxi|^2\leq2(1-\bxi\cdot\nne)$, term $\mathrm{IV}$ is controlled by~\eqref{energy b}, term $\mathrm{VI}$ is controlled by~\eqref{energy a}, and term $\mathrm{VIII}$ is controlled by~\eqref{energy c}. It remains to estimate terms $\mathrm{V}$ and  $\mathrm{VII}$ .  

For  term $\mathrm{V}$, use ~\eqref{pro} to decompose into
\begin{equation}
 \begin{aligned}
    \mathrm{V}=\underbrace{\int_\o\nabla\hh:\nne\otimes \nne|\nabla\pe|-\nabla\hh:\nne\otimes\nne\e
    |\Pi_\ue\nabla\ue|^2 dx}_{=:~\mathrm{V}^1}\\  +\underbrace{\sum_{i,j}\int_\o{(\nabla\hh)}_{i,j}\e\b((\p_i\ue-\Pi_\ue\p_i\ue)\cdot(\p_j\ue-\Pi_\ue\p_j\ue)\b)}_{=:~\mathrm{V}^2} dx.
\end{aligned}   
\end{equation}
The term $\mathrm{V}^2$ is controlled by~\eqref{energy c}.
For the term $\mathrm{V}^1$, apply \(|\nn_\e - \bxi|^2\leq2(1-\bxi\cdot\nne)\) and $|(\bxi\cdot\nabla)\hh|\leq C\dd_\Gamma\leq C\sqrt{1-\bxi\cdot\nne}$ (cf. \eqref{eq_calibration4}) to get
\begin{equation}
   \begin{aligned}
    |\mathrm{V}^1|=&\bigg|\int_\o\nabla\hh:\nne\otimes \nne\lef|\nabla\pe|-\e|\nabla\ue|^2\rig dx\bigg|\\
    =&\bigg|\int_\o\nabla\hh:\nne\otimes (\nne-\bxi)\lef|\nabla\pe|-\e|\nabla\ue|^2\rig dx\\
   & +\int_\o(\bxi\cdot\nabla)\hh\cdot\nne\lef|\nabla\pe|-\e|\nabla\ue|^2\rig dx\bigg|\\
   \leq&C\int_\o\sqrt{1-\bxi\cdot\nne}\left||\nabla\pe|-\e|\nabla\ue|^2\right|dx\\
   \leq&C\int_\o\sqrt{1-\bxi\cdot\nne}\left||\nabla\pe|-\e|\Pi_\ue
   \nabla\ue|^2\right|dx\\
   &+C\int_\o\sqrt{1-\bxi\cdot\nne}\left|\e|\nabla\ue|^2-\e|\Pi_\ue
   \nabla\ue|^2\right|dx.
\end{aligned} 
\end{equation}
The second integral in the last display can be bounded by~\eqref{energy c}. For the first integral, applying the  Cauchy-Schwarz inequality, we have
\begin{align*}
    &\int_\o\sqrt{1-\bxi\cdot\nne}\left||\nabla\pe|-\e|\Pi_\ue
   \nabla\ue|^2\right|dx\\
   =&\int_\o\sqrt{1-\bxi\cdot\nne}\e^{\frac{1}{2}}|\Pi_\ue
   \nabla\ue|\left| \e^{\frac{1}{2}}|\Pi_\ue
   \nabla\ue|-\e^{-\frac{1}{2}}|D\dd_F(\ue)|\right|dx\\
   \leq&\frac{1}{2}\int_\o(1-\bxi\cdot\nne)\e|\Pi_\ue
   \nabla\ue|^2 dx+\frac{1}{2} \int_\o\left| \e^{\frac{1}{2}}|\Pi_\ue
   \nabla\ue|-\e^{-\frac{1}{2}}|D\dd_F(\ue)|\right|^2dx,
\end{align*}
which is controlled by~\eqref{energy d} and~\eqref{energy e}. This completes the estimate for term $\mathrm{V}$.

For term $\mathrm{VII}$, use the inequality \( |a+b+c|^2 \leq 3|a|^2 + 3|b|^2 + 3|c|^2 \) to decompose as
\begin{align*}
     &\int_\o \Big|\e^{-\frac{1}{2}} \div \bxi  |D  \dd_F   (\uu_\e  )|\nn_\e  +\e^{\frac{1}{2}}  |\Pi_{\uu_\e  } \nabla \uu_\e  | \hh\Big|^2dx\\
     \leq&3\int_\o\Big|\div\bxi\lef\e^{-\frac{1}{2}}|D  \dd_F   (\uu_\e  )|-\e^{\frac{1}{2}}  |\Pi_{\uu_\e  } \nabla \uu_\e  | \rig\nne\Big|^2dx\\
     &+3\int_\o\Big|\div\bxi\e^{\frac{1}{2}}|\Pi_{\uu_\e  } \nabla \uu_\e  |(\nne-\bxi) \Big|^2dx\\
     &+3\int_\o\Big|(\div\bxi\bxi+\hh)\e^{\frac{1}{2}}|\Pi_{\uu_\e  } \nabla \uu_\e  | \Big|^2dx.
\end{align*}
The first integral on the right-hand side of the above inequality is controlled by~\eqref{energy d}. The second integral is controlled by~\eqref{energy e}, due to  $|\bxi-\nne|\leq2(1-\nne\cdot\bxi)$. The third integral can be treated using ~\eqref{energy f} combined with $\bxi\cdot\hh+\div\bxi=O(\dd_\Gamma)$ (cf.~\eqref{eq_calibrationEvolByMCF}).

This completes the proof of Proposition~\ref{proposition relative}, establishing the desired differential inequality for $\E$.
\end{proof}
\begin{corollary}\label{cor e estimate}
    There exists a constant $C_1=C_1(\Gamma)$ independent of $\e$ such that
    \begin{subequations}
    \begin{align}
        \sup_{t\in[0,T]}\E(t)&\leq C_1\e,\label{E estimate}
        \\
            \sup_{t\in[0,T]}\int_\o\b(|\nabla\pe|-\bxi\cdot\nabla\pe\b)dx&\leq C_1\e,\label{level e}\\
           \sup_{t\in[0,T]}\int_\o\b|\nabla\ue-\Pi_\ue\nabla\ue\b|^2dx&\leq C_1,\label{cor3.4 0.3}
           \\
        \int_0^T\int_\o\big|\p_t\ue+ (\hh\cdot\nabla)\ue\big|^2dxdt&\leq C_1 .\label{time estimate}
    \end{align}        
    \end{subequations}  
\end{corollary}
\begin{proof}
    We first apply the Cauchy-Schwarz inequality to obtain
    \begin{equation}\label{co3.4 1}
        \begin{aligned}
             &\e^2|\p_t\ue|^2-|\hh_\e|^2+\big|\hh_\e-\e\hh|\nabla\ue|\big|^2\\ =&\e^2|\p_t\ue|^2+\e^2|\hh|^2|\nabla\ue|^2+2\e^2(\hh\cdot\nabla)\ue\cdot\p_t\ue\\
        \geq& \e^2|\p_t\ue+(\hh\cdot\nabla)\ue|^2.
        \end{aligned}
    \end{equation}
    Combining the above inequality  with  the differential inequality in Proposition~\ref{proposition relative} and applying Grönwall’s inequality to the initial bound $\E(0)\leq C_0\e$, we deduce
\begin{equation}
    \begin{aligned}
        \sup_{t\in[0,T]}\E(t)&+\frac{\e}{2}\int_0^T\int_\o|\p_t\ue+(\hh\cdot\nabla)\ue|^2dxdt\\
        &\leq e^{C(\Gamma,T)}\E(0)\leq C_1\e,
    \end{aligned}
\end{equation}
which establishes~\eqref{E estimate} and~\eqref{time estimate}. The coercivity estimates~\eqref{energy c} and~\eqref{energy e}  then imply~\eqref{level e} and~\eqref{cor3.4 0.3}. This concludes the proof of Corollary~\ref{cor e estimate}.
\end{proof}

\begin{corollary}
    There exists a constant $C_1=C_1(\Gamma)>0$ independent of $\e$ such that \begin{subequations}
        \begin{align}
            \sup_{t\in[0,T]}\B(t)&\leq C_1\e\label{B ESTIMATE},\\
            \sup_{t\in[0,T]}\int_\o\big|\pe-c_F\chi_{\otp}\big|dx&\leq C_1\e^{\frac{1}{2}}.\label{unweighted}
        \end{align}
    \end{subequations}
\end{corollary}
\begin{proof}
The proof is divided into two steps. 

\textbf{Step $1$.  Proof of~\eqref{B ESTIMATE}.}\\
Define the difference function
\begin{equation}
    \omega_\e:=c_F\chi_{\o^+_t}-\pe=\frac{c_F}{2}\chi+\frac{c_F}{2}-\pe\quad\text{ where }\chi=\chi_\otp-\chi_\otn.
\end{equation}
By the generalized chain rule (cf.~\eqref{def general}),
\begin{equation}
    \p_t\pe=\big(\p_t\ue+(\hh\cdot\nabla)\ue\big)\cdot D\dd_F(\ue)-\hh\cdot\nabla\pe.
\end{equation}
Since $\var=0$ on $\Gamma$, we have $\frac{d}{dt} \int_\o \chi \var \,dx = 0$. A direct computation yields
\begin{align}
    \frac{d}{dt}\B
    =&-\int_\o \big(\p_t\ue
    +(\hh\cdot\nabla)\ue\big)\cdot D\dd_F(\ue)\var dx\\
    &+\int_\o\hh\cdot\nabla\pe\var dx
    +\int_\o\omega_\e\p_t\var dx.  
\end{align}
Observe that $(\p_i\chi)\var=0$, and thus $(\p_i\omega_\e)\var=-(\p_i\pe)\var$ in the sense of distribution. Using this and integrating by parts, we obtain
\begin{equation}
    \int_\o\hh\cdot\nabla\pe\var dx=-\int_\o\hh\cdot \nabla\omega_\e\var dx=\int_\o(\div\hh\var+\hh\cdot\var)\omega_\e dx.
\end{equation}
Using the Cauchy-Schwarz inequality, the generalized differential inequality $|D\df(\ue)|^2\leq2F(\ue)$ (cf.~\eqref{df a.e.}), and the transport equation of $\var$~\eqref{eq_weightEvol}, we conclude that
\begin{align}
     \frac{d}{dt}\B
     =&-\int_\o \big(\p_t\ue+(\hh\cdot\nabla)\ue\big)\cdot D\dd_F(\ue)\var dx\\
     &+\int_\o\div\hh\var\omega_\e dx
     +\int_\o(\p_t\var+\hh\cdot\var)\omega_\e dx\\
     \leq&\int_\o\frac{\e}{2}\big|\p_t\ue+(\hh\cdot\nabla)\ue\big|^2+\int_\o\frac{1}{\e}F(\ue)\var^2dx+C\B.
\end{align}
This inequality, combined with the estimate~\eqref{time estimate} and  the coercivity estimate~\eqref{energy f}   with the upper bound estimate of $\var$~\eqref{2.16d}, implies
\begin{equation}
    \frac{d}{dt}\B\leq C\e+C\B.\label{B estimate}
\end{equation}
Applying Grönwall’s inequality with the initial bound $\B(0)\leq C_0\e$ yields
\begin{equation}
    \sup_{t\in[0,T]}\B(t)\leq(\B(0)+CT\e )e^{CT}\leq C_1\e,
\end{equation}
thus proving~\eqref{B ESTIMATE}.

\textbf{Step $2$.  Proof of~\eqref{unweighted}.}\\
 We apply a slicing argument as in \cite{hensel2022convergence}.  Let $M\subset\rt$ be an embedded, compact, oriented $C^2$ curve, with unit normal $\nn_M$. Then, there exists a constant $r_M>0$ such that 
\begin{equation}
    \Psi_M:M\times(-r_M,r_M)\rightarrow\rt,\quad(x,s)\mapsto x+s\nn_M(x)\label{diffeo}
\end{equation}
is a $C^2$ diffeomorphism with bounds $|\nabla\Psi_M|, |\nabla\Psi_M^{-1}| \leq C_M$.

For $g\in L^\infty(\rt)$, the slicing inequality holds (cf.~\cite{fischer2020convergence}):
\begin{equation}
    \bigg|\int_{\Psi_M(M\times(-r_M,r_M))}|g|dx\bigg|^2\lesssim  \int_{\Psi_M(M\times(-r_M,r_M))}|g|\dist(\cdot,M)dx. \label{slicing 1}
\end{equation}

To proceed, we establish the following auxiliary claim.\\
\textit{Claim.} For any $g\in L^\infty(\o)$,
    \begin{equation}   \label{calim 3.49} \bigg|\int_\o|g|dx\bigg|^2\lesssim\int_\o|g|(\cdot)|\var|(\cdot,t)dx.
\end{equation}
\begin{proof}[Proof of the Claim]   
 Fix $t\in[0,T]$ and  define $r:=\min\{r_{\p\o},r_{\Gamma}\}$. Set
\begin{align}
    \o_{I}&:=\o\cap\Psi_{\Gamma_t}(\Gamma_t\times(-r,r)),\\
    \o_{O}&:=\{x\in\o:\dist(x,\Gamma_t)> r\}\cap\{x\in\o:\dist(x,\p\o)>r\},\\
    \o_{\p\o}&:=\o\setminus(\o_I\cup\o_O).
\end{align}
Since $\o=\o_{I}\cup \o_{O}\cup\o_{\p\o}$, by the inequality \( |a+b+c|^2 \leq 3|a|^2 + 3|b|^2 + 3|c|^2 \), \eqref{calim 3.49} decomposes as
\begin{equation}
    \bigg(\int_\o|g|dx\bigg)^2\leq 3\bigg(\int_{\o_{I}}|g|dx\bigg)^2+3\bigg(\int_{\o_{O}}|g|dx\bigg)^2+3\bigg(\int_{\o_{\p\o}}|g|dx\bigg)^2.
\end{equation}
For the first term in the right-hand side, due to $\dist(\cdot,\Gamma_t)\leq C\dist(\cdot,\p\o)$, we apply the slicing estimate~\eqref{slicing 1} with $M=\Gamma_t$ and the lower bound~\eqref{eq_weightCoercivity} to estimate it. For the second term, because of the uniform lower bound $r$ of $\dist(\cdot,\Gamma_t)$ and $\dist(\cdot,\p\o)$ in $\o_O$, we can directly apply~\eqref{eq_weightCoercivity}. For the third term, noting that $\o_{\p\o}\subset\Psi_{\p\o}(\p\o\times(-r,r))$ and $\dist(\cdot,\p\o)\leq C\dist(\cdot,\Gamma_t)$ in $\o_{\p\o}$ (due to the boundary contact angle condition $|\no\cdot\nn_{\Gamma_t}|=\cos\alpha>0$), we also control it by~\eqref{eq_weightCoercivity}.
\end{proof}
Applying this claim to $g=\pe-c_F\chi_{\otp}$ and using the estimate~\eqref{B ESTIMATE}, we complete the proof of~\eqref{unweighted}.
\end{proof}

\section{Proof of Theorem~\ref{theorem 2}}\label{proof the 2}

\subsection{Weak convergence techniques}
Combining Corollary~\ref{cor e estimate} and the coercivity estimate~\eqref{energy f}, we obtain the  following results.
\begin{corollary}\label{coro energy}
 There exists a constant $C=C(\Gamma)$ independent of $\e,\delta>0$ such that the following estimates hold:
\begin{subequations}
\begin{align}
        \|\ue\|_{L^\infty(\o\times(0,T))}&\leq R_N,\label{boundness}\\
        \sup_{t\in[0,T]}\int_{\o_t^\pm\setminus B_\delta(\Gamma_t)}\lef\frac{1}{2}|\nabla\ue|^2+\frac{1}{\e^2}F(\ue)\rig dx&\leq C\delta^{-2},\label{estimate b}\\
        \int_0^T\int_{\o_t^\pm\setminus B_\delta(\Gamma_t)}|\p_t\ue|^2dxdt&\leq C\delta^{-2}.  \label{coro energy 2}
\end{align}
\end{subequations}
\end{corollary}
\begin{proof}
The uniform $L^\infty$-bound estimate~\eqref{boundness} is directly derived from Lemma~\ref{ue regularity}.  
The estimate~\eqref{estimate b} is obtained by combining the  energy bound~\eqref{E estimate} 
with the coercivity estimate~\eqref{energy f}.  Finally,~\eqref{coro energy 2} is deduced from the time–derivative bound
\eqref{time estimate} along with the same coercivity estimate~\eqref{energy f}.
\end{proof}
In particular, for each fixed $\delta>0$, the sequence $\{\ue\}$ is uniformly bounded in
\begin{equation}
     L^\infty\Big(0,T;L^\infty(\Omega)\cap H^1\b(\o_t^\pm\setminus B_\delta(\Gamma_t)\b)\Big)\bigcap H^1\Big(0,T;L^2\b(\o_t^\pm\setminus B_\delta(\Gamma_t)\b)\Big).
\end{equation}

\begin{corollary}\label{weak con}
    For every sequence $\e_k\downarrow0$ there exist a subsequence (not  relabel) and limits
    $$\uu^\pm\in L^\infty\big(0,T;L^\infty(\o)\cap H^1_{\mathrm{loc}}(\o_t^{\pm},N^\pm)\big), \qquad \p_t\uu^\pm\in L^2\big(0,T;L^2_{\mathrm{loc}}(\o_t^\pm)\big),$$ such that
    \begin{subequations}
        \begin{align}    
 \p_t\uu_{\e_k}&\xrightharpoonup{k\to\infty}\p_t\uu^\pm &&\text{ weakly in } L^2\big(0,T;L^{2}_{\mathrm{loc}}(\o_t^{\pm})\big),\label{limit 1}\\        \nabla\uu_{\e_k}&\xrightharpoonup{k\to\infty}\nabla \uu^\pm
       &&\text{ weakly-star in }  L^\infty\big(0,T;L^2_{\mathrm{loc}}(\o_t^{\pm})\big),\label{limit 2}\\
        \uu_{\e_k}&\xrightarrow{k\to\infty}\uu^\pm
        &&\text{ strongly in } C\big([0,T];L^2_{\mathrm{loc}}(\o_t^{\pm})\big),\label{limit 3}\\
       \uu_{\e_k}&\xrightarrow{k\to\infty}\uu^\pm
        &&\text{ a.e. in } \bigcup_{t\in[0,T]}\o_t^\pm\times\{t\}.\label{limit 4}
    \end{align}
    \end{subequations}   
\end{corollary}
\begin{proof}
     Fix $\delta>0$. By using the weak-star and weak compactness and the Aubin-Lions lemma, there exists a subsequence  $\e_k=\e_k(\delta)$ such that
     \begin{subequations}
         \begin{align}
          \uu_{\e_k}&\xrightharpoonup{k\to\infty}\uu
        &&\text{ weakly-star in } L^{\infty}\big(0,T;L^{\infty}(\o)\big)\label{cor4.2 3},\\
\p_t\uu_{\e_k}&\xrightharpoonup{k\to\infty}\p_t\uu_\delta^\pm &&\text{ weakly in } L^2\big(0,T;L^2(\o_t^\pm\setminus B_\delta(\Gamma_t))\big),\label{cor4.2 1}\\
        \nabla\uu_{\e_k}&\xrightharpoonup{k\to\infty}\nabla \uu_\delta^\pm
       &&\text{ weakly-star in }  L^\infty\big(0,T;L^2(\o_t^\pm\setminus B_\delta(\Gamma_t))\big),\label{cor4.2 2}\\
         \uu_{\e_k}&\xrightarrow{k\to\infty}\uu^\pm_\delta
        &&\text{ strongly in } C\big([0,T];L^2(\o_t^\pm\setminus B_\delta(\Gamma_t))\big),\\
        \uu_{\e_k}&\xrightarrow{k\to\infty}\uu^\pm_\delta
        &&\text{ a.e. in }\bigcup_{t\in(0,T)}\big(\o_t^\pm\setminus B_\delta(\Gamma_t)\big)\times\{t\}.
         \end{align}
     \end{subequations}
      Moreover, $\uu=\uu_\delta^\pm$ a.e. in $\bigcup_{t\in(0,T)}\big(\o_t^\pm\setminus B_\delta(\Gamma_t)\big)\times\{t\}$. To obtain convergence independent of $\delta$, we use a standard diagonal argument to extract a common subsequence that satisfies~\eqref{limit 1}-\eqref{limit 4}.  
      
     We proceed to determine $\uu^\pm$. Fix $t\in[0,T]$ and let $K\subset\o_t^\pm$ be a compact set. Then, by Fatou's lemma combined with  the pointwise convergence~\eqref{limit 4} and the estimate~\eqref{estimate b}, we have
     \begin{equation}
         \int_KF(\uu)dx\leq\liminf_{k\to\infty}\int_K F(\uek)dx=0,
     \end{equation}
     which implies that $\uu\in N^\pm$ a.e. in $\o\times[0,T]$. On the other hand, since $\df$ is Lipschitz continuous (cf. Lemma~\ref{df}), from the strong convergence~\eqref{limit 3}, we obtain
     \begin{equation}  \label{df regularity}      \df(\uek)\xrightarrow{k\to\infty}\df(\uu^\pm)\quad\text{ strongly in } C\big([0,T];L_{\mathrm{loc}}^2(\o_t^\pm)\big).
     \end{equation}
Combining this with the bulk estimate~\eqref{unweighted} and the property of $\df$~\eqref{df pro2}, we conclude that  $\uu^\pm$ maps $\o_t^\pm$ into $N^\pm$ respectively.  
\end{proof}

\subsection{\texorpdfstring{Improve the regularity of $\uu^\pm$ in $\o_t^\pm$}{}} \label{improve the regularity}
In order to ensure the well-definedness of the limits  \(\uu^\pm\) on the interface \(\g\), we proceed to upgrade the regularity of $\uu^\pm$ from $\uu^\pm\in L^{\infty}\big(0,T;H^1_{\mathrm{loc}}(\Omega_t^\pm;N^\pm)\big)$ to $\uu^\pm\in L^{\infty}\big(0,T;H^1(\Omega_t^\pm;N^\pm)\big)$. 

 Since $\df$ is Lipschitz continuous, we fix a small constant $\delta_0\in(0,\frac{\cf}{4})$ such that
\begin{equation}\label{de 0}
    \begin{aligned}
     \{\uu\in\rk:\df(\uu)<4\delta_0\}\subset B_{\delta_N}(N^-),\\
    \{\uu\in\rk:\df(\uu)>c_F-4\delta_0\}\subset B_{\delta_N}(N^+), 
    \end{aligned}
\end{equation}
where $\delta_N$ is given in~\eqref{delta_0}.

The following proposition is used to derive the level set estimates in Proposition~\ref{level set estimate}.

\begin{proposition}\label{level set}
      There exists a constant $C=C(\Gamma)$  independent of $(\e,t,\delta)$ such that the following statements hold:     
      \begin{itemize}
             
      \item[(1)] For any $\delta\in(0,4\delta_0)$,
\begin{equation}
    \big|\{\psi_{\e}
    >\delta\}\Delta\otp\big|+\big|\{\pe<\cf-\delta\}\Delta\otn\big|\leq C\delta^{-1}\e^\frac{1}{2}.\label{bulk error}
\end{equation}
\item[(2)]  For  a.e. $\delta\in(0,4\delta_0)$, the sets $\{\psi_{\e}
    >\delta\}$ and $\{\pe<\cf-\delta\} $  are sets of finite perimeter in $\o$ and $\rn$, and
      \begin{align}
          \hn\Big( \big(\p^*\{\psi_{\e}
    >\delta\}\cap\p\o\big)\Delta\big(\p\otp\cap\p\o\big)\Big)&\leq C\delta^{-1}\e^\frac{1}{2},\label{pr4.3 1.1}\\
    \hn\Big(\big(\p^*\{\pe<\cf-\delta\}\cap\p\o\big)\Delta\big(\p\otn\cap\p\o\big)\Big)&\leq C\delta^{-1}\e^\frac{1}{2}.\label{pr4.3 1.2}
      \end{align}
      \end{itemize}
 

\end{proposition}


\begin{proof}
The proof is divided into two steps.

\textbf{Step $1$.  Proof of~\eqref{bulk error}.}\\
Using the bulk estimate~\eqref{unweighted} and the Chebyshev inequality, we have
\begin{equation}
    \begin{aligned}
         &\big|\{\psi_{\e}
    >\delta\}\Delta\otp\big|\\=&\big|\{\psi_{\e}
    >\delta\}\cap\otn\big|+\big|\{\pe\leq\delta\}\cap\otp\big|\\
        \leq&\delta^{-1}\int_{\otn\cap\{\pe>\delta\}}\pe dx+(c_F-\delta)^{-1}\int_{\otp\cap\{\pe\leq\delta\}}\big(c_F-\pe\big) dx\\
        \leq&C\delta^{-1}\e^{\frac{1}{2}}.
    \end{aligned}
\end{equation}
We can estimate $\big|\{\pe<\cf-\delta\}\Delta\otn\big| $ analogously to obtain the estimate~\eqref{bulk error}.

\textbf{Step $2$. Proofs of~\eqref{pr4.3 1.1} and~\eqref{pr4.3 1.2}.}\\
By the bulk estimate~\eqref{unweighted} and the Sobolev trace theorem, we have
\begin{equation}\label{boundary error}
\begin{aligned}
    \sup_{t\in[0,T]} &\int_{\p\o
    } (c_F\chi_{\o_t^+}-\psi_\e) d\hn\leq C\e^{\frac{1}{2}}.
\end{aligned}
\end{equation}
 From the co-area formula for $BV$ functions, for a.e. $\delta\in(0,\delta_0)$, the set $\{\pe>\delta\}$ is of finite perimeter in $\o$, and by Lemma~\ref{finite c}, also in $\R^n$. For $\hn$-a.e. $x\in \p^*\{\pe>\delta\}\cap\p\o$, we use the local trace property (cf. \cite[Theorem 5.7]{evans2018measure}) and the density property of the reduced boundary (cf. \cite[Corollary 15.8]{maggi2012sets}) to obtain
\begin{equation}
\begin{aligned}
     \pe(x,t)=&\lim_{r\downarrow0}\frac{\int_{B_r(x)\cap\o}\pe(y,t)dy}{|B_r(x)\cap\o|}\\
     \geq&\limsup_{r\downarrow0}\frac{\int_{B_r(x)\cap\{\pe>\delta\}}\pe(y,t)dy}{|B_r(x)\cap\{\pe>\delta\}|}\frac{|B_r(x)\cap\{\pe>\delta\}|}{|B_r(x)\cap\o|}\geq\delta.
\end{aligned}\label{pr4.3 3}
\end{equation}
This implies that 
\begin{equation}\label{pr4.3 4}
    \pe(x,t)\geq\delta\quad\text{ for }\hn\text{ -a.e. }x\in \p^*\{\pe>\delta\}\cap\p\o.
\end{equation}
By a similar argument, we deduce
\begin{equation}\label{pr4.3 5}
    \pe(x,t)\leq\delta\quad\text{ for }\hn\text{ -a.e. }x\in \p^*\{\pe\leq\delta\}\cap\p\o.
\end{equation}
Combining~\eqref{boundary error},~\eqref{pr4.3 4}, and~\eqref{pr4.3 5} with  Lemma~\ref{finite c}, the  Chebyshev inequality yields
\begin{equation}
    \begin{aligned}
        &\hn\big((\p^*\{\psi_{\e}
    >\delta\}\cap\p\o)\Delta(\p\otp\cap\p\o)\big)\\
    =&\hn(\p^*\{\pe>\delta\}\cap\p\otn\cap\p\o)+\hn(\p^*\{\pe\leq\delta\}\cap\otp\cap\p\o)\\
    \leq&\delta^{-1}\int_{\p \otn\cap\p\o}\pe d\hn+(c_F-\delta)^{-1}\int_{\p\otp\cap\p\o}(c_F-\pe) d\hn\\
    \leq& C\delta^{-1}\e^{\frac{1}{2}}.
    \end{aligned}
\end{equation}
The estimate~\eqref{pr4.3 1.2} can be derived in the same manner.

This concludes the
proof of Proposition~\ref{level set}.
\end{proof}
The following result estimates the discrepancy between the interior reduced boundary of upper/lower level sets of \(\psi_\e(\cdot,t)\) and the interface $\Gamma_t$.
\begin{proposition}\label{level set estimate}
For any fixed triple $(\e,\delta,t)$ with $\delta<\delta_0$, there exist $a^\pm_{\e,\delta}(t)\in(\delta,4\delta)$ and a constant $C=C(\Gamma)>0$ independent of $(\e,t,\delta)$ such that 
    
    \begin{align}
         \left|\hn\lef \p^*\{\psi_{\e}< a^-_{\e,\delta}(t)\}\cap\o\rig -\hn(\Gamma_t)\right|&\leq C\delta^{-2}\e^\frac{1}{2},\label{pr4.4 1.1}\\
          \left|\hn\lef \p^*\{\psi_{\e}>c_F- a^+_{\e,\delta}(t)\}\cap\o\rig -\hn(\Gamma_t)\right|&\leq C\delta^{-2}\e^\frac{1}{2}.\label{pr4.4 1.2}
    \end{align}
     
\end{proposition}
\begin{proof}
To begin, we choose $\pe(\cdot,t)=\pe^*(\cdot,t)\in H^1(\o)$ (see \eqref{df regularity} for its regularity) as its precise representative (as defined in~\eqref{preciser resprestive}) and recall the regularity of $\bxi$ from ~\eqref{eq_regularityXi}. By Lemmas~\ref{finite c} and~\ref{precise}, for any fixed triple $(\e,\delta,t)$, we may assume there exist $\delta^-_\e=\delta_\e^-(t,\delta)\in(\delta,2\delta)$ and  $\delta^+_\e=\delta_\e^+(t,\delta)\in(3\delta,4\delta)$  such that the disjoint sets 
\begin{align*}
      \o_{t}^{\e,-}&:=\{\psi_{\e}(x,t)<\delta^-_\e\},\\
     \o_{t}^{\e}&:=\{\delta^-_\e\leq\psi_{\e}(x,t)\leq\delta^+_\e\}, \\
     \o_{t}^{\e,+}&:=\{\psi_{\e}(x,t)>\delta^+_\e\},   
\end{align*}
have finite perimeter in $\o$ and $\rn$, and satisfy
\begin{equation}
\begin{aligned}
\p^*\big(\oten\cup\otep)\cap\o\overset{\hn}{=}(\p^*\oten\cap\o)\bigcup(\p^*\otep\cap\o).\label{reduced eq}
\end{aligned}
\end{equation}
Additionally, we have the following boundary conditions
\begin{subequations}\label{finite}
    \begin{align}
         \pe(x,t)&=\delta^-_\e\quad&&\hn\text{-a.e. }x\in\p^*  \o_{t}^{\e,-}\cap\o,\label{finite 1}\\
        \pe(x,t)&=\delta^+_\e\quad&&\hn\text{-a.e. }x\in\p^*\o_{t}^{\e,+}\cap\o,\label{finite 2}\\
        \bnu^\pm&=-\bnu\quad&&\hn\text{-a.e. }x\in\p^*\o_t^{\e,\pm}\cap\p\o_t^\e\cap\o,\label{finite 3}\\
        \bnu^\pm&=-\no \quad&&\hn\text{-a.e. }x\in\p^*\o_{t}^{\e,\pm}\cap\p\o,\label{finite 4}
    \end{align}
\end{subequations}
where $\bnu$ denotes the measure-theoretic outer unit normal to $\p^*\oted$, and $\bnu^\pm$ denotes the  measure-theoretic outer unit normal to $\p^*\o_t^{\e,\pm}$.

 By the generalized Gauss-Green formula (Lemma~\ref{Generalized Gauss-Green formula}), the co-area formula  and estimate~\eqref{level e}, we derive
 \begin{equation}\label{pr4.4 1}
     \begin{aligned}
          C_1{\e}\geq&\int_{\o_t^{{\e}}}\b(|\nabla\psi_{\e}|-\bxi\cdot\nabla\psi_{\e} \b)dx\\
    \geq&\int_{{\delta^-_\e}}^{\delta^+_\e}\hn\lef \p^*\{\psi_{\e}<s\}\cap\o\rig ds\\
    &-\int_{\p^*\o_t^{\e}\cap\o}\bxi\cdot\bnu \psi _{\e} d\hn+\int_{\p^*\o_t^{\e}\cap\p\o}\bxi\cdot\no \psi _{\e} d\hn+\int_{\o_t^{{\e}}}\div\bxi\psi_{\e} dx.
     \end{aligned}
 \end{equation}

We proceed in two steps.

 \textbf{Step $1$. Estimates for $ \int_{\p^*\o_t^{\e}\cap\p\o}\bxi\cdot\no \psi _{\e} d\hn$ and $\int_{\o_t^{{\e}}}\div\bxi\psi_{\e} dx$.}\\
 We first observe that $\pe\in[\delta,\cf-\delta]$ in $\oted$. By applying the Chebyshev inequality and the bulk estimate~\eqref{unweighted}, we have
\begin{equation}\label{pr4.4 2}
    \bigg|\int_{\o_t^{{\e}}}\div\bxi\psi_{\e} dx\bigg|\leq C|\oted|\leq C\delta^{-1}\int_\o(c_F\chi-\pe)dx\leq C\delta^{-1}\e^{\frac{1}{2}}.
\end{equation}
By the trace estimate as in~\eqref{pr4.3 3}, $\pe\in[\delta,c_F-\delta]$ on $\p^*\oted\cap\p\o$. Then, by using the Chebyshev inequality and the boundary estimate~\eqref{boundary error}, we obtain 
    \begin{equation}  \label{pr4.4 7}
    \begin{aligned}
          \bigg|\int_{\p^*\o_t^{\e}\cap\p\o}\bxi\cdot\no \psi _{\e} d\hn\bigg|\leq& C\hn(\p^*\oted\cap\p\o)\\
          \leq&C\delta^{-1}\int_{\p\o}(c_F\chi-\pe)d\hn\leq C\delta^{-1}\e^{\frac{1}{2}}.
    \end{aligned}
\end{equation}
Combining~\eqref{pr4.4 1}-\eqref{pr4.4 7}, we obtain
\begin{align}
    \left|\int_{{\delta^-_\e}}^{{\delta^+_\e}}\hn\lef \p^*\{\psi_{\e}<s\}\cap\o\rig ds-\int_{\p^*\o_t^{{\e}}\cap\o}\bxi\cdot\bnu \pe d\hn\right|
   \leq C\delta^{-1}\e^{\frac{1}{2}}. \label{error term 1}
\end{align}

\textbf{Step $2$. Estimate for $\int_{\p^*\o_t^{{\e}}\cap\o}\bxi\cdot\bnu \pe d\hn$.}\\
Using the conditions ~\eqref{finite 1} and~\eqref{finite 2} combined with~\eqref{reduced eq}, the boundary integral decomposes as 
\begin{equation}\label{pr4.4 3}
    \begin{aligned}
          &\int_{\p^*\o_t^{{\e}}\cap\o}\bxi\cdot\bnu \psi_{\e} d\hn\\
         =&\delta^+_\e\int_{\p^*\o_{t}^{\e,+}\cap\o} \bxi\cdot\bnu d\hn+\delta^-_\e\int_{\p^*\o_{t}^{\e,-}\cap\o} \bxi\cdot\bnu d\hn.
    \end{aligned}
\end{equation}
By applying the generalized Gauss-Green formula and using $\bnu^\pm=-\bnu$ on $\p^*\o_t^{\e,\pm}\cap\p\o_t^\e\cap\o$ (cf.~\eqref{finite 3}), we obtain
 \begin{align}
  \delta^+_\e \int_{\o_{t}^{\e,+}}&\div\bxi dx=\nonumber\\-\delta^+_\e&\int_{\p^*\o_{t}^{\e,+}\cap\o}\bxi\cdot\bnu d\hn-\delta^+_\e\int_{\p^*\o_{t}^{\e,+}\cap\p\o}\bxi\cdot\no d\hn,\label{e inner 1}\\
    \delta^-_\e\int_{\o_{t}^{\e,-}}&\div\bxi dx=\nonumber\\-{\delta_\e^-}&\int_{\p^*\o_{t}^{\e,-}\cap\o}\bxi\cdot\bnu d\hn-{\delta_\e^-}\int_{\p^*\o_{t}^{\e,-}\cap\p\o}\bxi\cdot\no d\hn.\label{e inner 2}
 \end{align}
Substituting~\eqref{e inner 1}-\eqref{e inner 2} into~\eqref{pr4.4 3}, we find
\begin{equation}\label{pr4.4 5}
\begin{aligned}
    &\int_{\p^*\o_t^{{\e}}\cap\o}\bxi\cdot\bnu \psi_{\e} d\hn\\
    =&-\delta^+_\e \int_{\o_{t}^{\e,+}}\div\bxi dx-{\delta^-_\e}\int_{\o_{t}^{\e,-}}\div\bxi dx\\   &-\delta^+_\e\int_{\p^*\o_{t}^{\e,+}\cap\p\o}\bxi\cdot\no d\hn-{\delta^-_\e}\int_{\p^*\o_{t}^{\e,-}\cap\p\o}\bxi\cdot\no d\hn.
\end{aligned}
\end{equation}
On the other hand, applying the standard Gauss-Green formula to $\div\bxi$ in $\o_t^\pm$, we have
\begin{equation}\label{pr4.4 6}
\begin{aligned}
    -({\delta_\e^+}-\delta_\e^-)\hn(\Gamma_{t})=&\delta_\e^+\int_{\o^+_t}\div\bxi dx+\delta_\e^-\int_{\o^-_t}\div\bxi dx\\
    &+\delta_\e^+\int_{\p\o^+_t\cap\p\o}\bxi\cdot\no d\hn+\delta_\e^-\int_{\p\o^-_t\cap\p\o}\bxi\cdot\no d\hn.
\end{aligned}  
\end{equation}
Comparing~\eqref{pr4.4 5} with~\eqref{pr4.4 6}, we arrive at
\begin{equation}\label{pr4.4 4}
\begin{aligned}
     &\left|\int_{\p^*\o_t^{{\e}}\cap
     \o}\bxi\cdot\bnu \psi_{\e} d\hn- {(\delta_\e^+-\delta_\e^-)}\hn(\Gamma_{t})\right|\\
     \leq&\underbrace{\delta_\e^+\left|\int_{\o_{t}^{\e,+}}\div\bxi dx-\int_{\o^+_t}\div\bxi dx\right|}_{=:~i}+\underbrace{\delta_\e^-\left|\int_{\o_{t}^{\e,-}}\div\bxi dx-\int_{\o^-_t}\div\bxi dx\right|}_{=:~ii}\\
     &+\underbrace{\delta_\e^+\left|\int_{\p\o^+_t\cap\p\o}\bxi\cdot\no d\hn-\int_{\p^*\o_{t}^{\e,+}\cap\p\o}\bxi\cdot\no d\hn\right|}_{=:~iii}\\
     &+\underbrace{\delta_\e^-\left|\int_{\p\o^-_t\cap\p\o}\bxi\cdot\no d\hn-\int_{\p^*\o_{t}^{\e,-}\cap\p\o}\bxi\cdot\no d\hn\right|}_{=:~iv}.
\end{aligned}
\end{equation}

Each of the error terms is estimated as follows.

For terms $i$ and $ii$, by the bulk error estimate~\eqref{bulk error} in Proposition~\ref{level set estimate}, 
    \begin{equation}
        i+ii\leq C\lef|\otp\Delta\o_{t}^{\e,+}|+|\otn\Delta\o_{t}^{\e,-}|\rig\leq C\delta^{-1}\e^{\frac{1}{2}}.
    \end{equation}
    
For terms $iii$ and $iv$, by the boundary error estimates ~\eqref{pr4.3 1.1}-\eqref{pr4.3 1.2} in Proposition~\ref{level set estimate}, 
    \begin{equation}
    \begin{aligned}
       iii+iv \leq &C\hn\Big( \big(\p^*\otep\cap\p\o\big)\Delta\big(\p\otp\cap\p\o\big)\Big)\\
    &+C\hn\Big(\big(\p^*\oten\cap\p\o\big)\Delta\big(\p\otn\cap\p\o\big)\Big)\\
    \leq&C\delta^{-1}\e^{\frac{1}{2}}. 
    \end{aligned}         
    \end{equation}

Combining the above estimates and dividing both sides of~\eqref{error term 1} by $(\delta_\e^+-\delta^-_\e)\geq\delta$, we obtain
\begin{align}
    \left|\frac{1}{{\delta_\e^+-\delta_\e^-}}\int_{\delta_\e^-}^{\delta_\e^+}\hn\lef \p^*\{\psi_{\e}<s\}\cap\o\rig  ds-\hn(\Gamma_t)\right|\leq C\delta^{-2}\e^{\frac{1}{2}}.
\end{align}

Finally, by Fubini's theorem, we conclude that there exists $a_{\e,\delta}^-(t) \in(\delta_\e^-,\delta_\e^+)\subset (\delta, 4\delta)$ satisfying the first estimate~\eqref{pr4.4 1.1}.  The second estimate~\eqref{pr4.4 1.2} is obtained by a symmetric argument, considering the  upper/lower level sets with $\delta_\e^-\in(c_F-4\delta,c_F-3\delta)$ and $\delta_{\e}^+\in(c_F-2\delta,c_F-\delta)$. Specifically, the estimates~\eqref{pr4.4 2} and~\eqref{pr4.4 7} hold relying only on the a priori bound $\pe\in[\delta,c_F-\delta]$ on the relevant sets.

This concludes the proof of Proposition~\ref{level set estimate}.

\end{proof} 
The next proposition establishes a relationship between  \(\nabla\uu\) and its tangential component on $N$.
 \begin{proposition}\label{pr core}
    For $\uu\in H^1(\o;\rk)$, there exists a constant $C=C(N)$ such that for $0<\delta<\delta_N$, where $\delta_N$ is given in~\eqref{delta_0}, and for a.e. $x\in\{x\in\o:\uu(x)\in B_{\delta}(N)\}$, we have
     \begin{align}
          |\nabla\uu-\Pi_\uu\nabla\uu|^2\geq(1-C\delta)\b|\nabla\b(P_N(\uu)\b)\b|^2,\label{core 1}
     \end{align}  
     where $P_N$ is the nearest point projection and $\Pi_\uu$ is defined in~\eqref{Pi}.
 \end{proposition}
  \begin{proof}
      To begin, fix any $ B_\dprr(y_0)\subset B_{\delta}(N)$ with some constant $\delta^{\prime\prime}>0$. It  suffices to verify~\eqref{core 1} for a.e. $x\in\{\uu\in B_\dprr(y_0)\}$. We choose $\delta^{\prime}>\dprr$ such that $B_\dprr(y_0)\subset B_{\delta^{\prime}}(y_0)\subset B_{\delta}(N)$  and consider  the radial truncation $ T_{\delta^{\prime},y_0}\in Lip(\rk;\rk)$  (see Lemma~\ref{radial trun} for its properties) defined by
      \begin{equation}
        T_{\delta^{\prime},y_0}(y)=\begin{cases}
            y&\text{ if }|y-y_0|\leq \delta^{\prime},\\
            y_0+ \delta^{\prime}\frac{y-y_0}{|y-y_0|}&\text{ if } |y-y_0|>\delta^{\prime}.
        \end{cases}
    \end{equation}
      Since $T_{\delta^{\prime},y_0}$ is Lipschitz, $\uu \in H^1(\o)$, and the Jacobian $D T_{\dpr,y_0}$ is the identity matrix where $\uu\in B_{\dpr}(y_0)$, then the generalized chain rule for weak derivatives (cf.~\eqref{general}-\eqref{def general}) gives 
\begin{align*}
     \nabla \big(T_{\delta^{\prime},y_0}(\uu)\big)=\nabla\uu \quad\text{ for a.e. }x\in\{\uu\in B_\dprr(y_0)\}.
\end{align*}
Therefore, we may replace $\uu$ with $T_{\delta^{\prime}}(\uu)$ without changing the gradient on the set of interest, reducing the proof to the case 
      \begin{equation}
          \uu(x)\in B_\dpr(y_0)\quad\text{ for a.e. }x\in\o. 
      \end{equation}
Without loss of generality, we assume $B_\dpr(y_0)\subset B_{\delta_N}(N^+)$.  Let $l\in\mathbb{N}^+$ be the co-dimension of $N^+$ and  $\b\{\nj\b(P_N(\cdot)\b)\b\}_{j=1}^l\subset C^1\b(  B_\dpr(y_0)\b)$ be a local orthonormal frame of the normal space in $P_N\b( B_\dpr(y_0)\b)$. Let  $$\{\dd_j(\cdot)\}_{j=1}^l\subset C^1\b(  B_\dpr(y_0)\b)$$ be the coordinates such that
    \begin{align}
          \uu=P_N\uu+\dd_N(\uu)\bnu_N(\uu)=P_N\uu+\sum_{j=1}^l\dd_j(\uu)\bnu_j(P_N\uu), \label{descomp}
      \end{align}
      where $\bnu_N$ is defined in~\eqref{bnu}, and
      \begin{equation}
          \dd_N^2(\uu)=\sum_{j=1}^l\dd_j^2(\uu) \text{ and }D\dd_j(\uu)=\bnu_j(P_N\uu) \quad\forall x\in \o,\label{core 1.1}
      \end{equation}
      where $\dd_N$ is the distance to $N$ (cf.~\eqref{distance to N}).
     By the chain rule (cf. \cite[Theorem 4.4]{evans2018measure}),  for a.e. $x\in\o$, we have 
      \begin{equation}\label{pr4.5 8}
           \begin{aligned}          \p_i\uu&=\p_i(P_N\uu)+\sum_{j=1}^l\Big(\p_i\big(\dd_j(\uu)\big)\bnu_j(P_N\uu)+\dd_j(\uu)\p_i\big(\bnu_j(P_N\uu)\big)\Big).\\       
      \end{aligned}  
      \end{equation}
   The subsequent computations hold for a.e. $x \in \o$.
      Note that,
      \begin{equation}\label{pr4.5 9}
          \p_i(P_N\uu)\cdot\bnu_j(P_N\uu)=0 \text{ and }\p_i\big(\bnu_j(P_N\uu)\big)\cdot\bnu_j(P_N\uu)=0.
      \end{equation}
      Using these relations, $\p_i\big(\dd_j(\uu)\big) = \bnu_j \cdot \p_i\uu$ from~\eqref{core 1.1} and the orthogonality $\bnu_j \cdot \bnu_k = \delta_{j,k}$, then~\eqref{pr4.5 8} gives
      \begin{equation}
      \begin{aligned}
          |\p_i\uu|^2=&|\p_i(P_N\uu)|^2+\bigg|\sum_{j=1}^l\p_i\big(\dd_j(\uu)\big)\bnu_j(P_N\uu)\bigg|^2\\
          &+2\p_i(P_N\uu)\cdot\sum_{j=1}^l\dd_j(\uu)\p_i\big(\bnu_j(P_N\uu)\big)\\
          =&|\p_i(P_N\uu)|^2+\sum_{j=1}^l|\bnu_j(P_N\uu)\cdot\p_i\uu|^2\\
          &+2\sum_{j=1}^l \dd_j(\uu)A^j(P_N\uu)\big(\p_i(P_N\uu),\p_i(P_N\uu)\big),
      \end{aligned}\label{core 2}
      \end{equation}
where $A^j(P_N\uu)(\vv,\ww):=D\big(\bnu_j(P_N)\big)(P_N\uu)\vv\cdot\ww$ is the j-th component of the second fundamental form of $N$ at $P_N\uu$. 
Since the second  fundamental form of $N$ is bounded, i.e., $|A^j| \leq C(N)$, where  $C(N)$ depends only on the geometry of $N$, and $\dd_j(\uu) \leq \delta$ due to $\uu \in B_\delta(N)$, we have
\begin{equation}
     \bigg |2\sum_{j=1}^l \dd_j(\uu)A^j\b(\p_i(P_N\uu), \p_i(P_N\uu)\b)\bigg|
    \leq  C(N) \delta |\p_i(P_N\uu)|^2.
\end{equation}  
Thus, the identity~\eqref{core 2} gives
\begin{equation}\label{pr4.5 1}
    |\p_i\uu|^2\geq\b(1-C(N)\delta\b)|\p_i(P_N\uu)|^2+\sum_{j=1}^l|\bnu_j(P_N\uu)\cdot\p_i\uu|^2.
\end{equation}

On the other hand, if $\uu(x)\in N$, then $D\df(\uu)(x)=\mathbf{0}$ and $\Pi_{\uu}\p_i\uu=\mathbf{0}$, so the relation~\eqref{core 1} follows  immediately from~\eqref{core 2}. If instead $\uu(x)\notin N$, by the property of $\df$~\eqref{df pro3} and the relations in~\eqref{core 1.1}, the Cauchy-Schwartz inequality gives
\begin{equation}
\begin{aligned}
     |\Pi_{\uu}\p_i\uu|^2\dd_N^2(\uu)&=|\p_i\uu\cdot\bnu_N(\uu)\dd_N(\uu)|^2\\
    &=\bigg|\sum_{j=1}^l\dd_j(\uu)\p_i\uu\cdot\bnu_j(P_N\uu) \bigg|^2\\
    &\leq \sum_{j=1}^l \dd_j^2(\uu)\sum_{j=1}^l|\p_i\uu\cdot\bnu_j(P_N\uu)|^2\\
    &=\dd_N^2(\uu)\sum_{j=1}^l|\p_i\uu\cdot\bnu_j(P_N\uu)|^2,
\end{aligned}\label{core 3}
\end{equation}
which implies
\begin{equation}\label{core 4}
    |\Pi_{\uu}\p_i\uu|^2=\sum_{j=1}^l|\p_i\uu\cdot\bnu_j(P_N\uu)|^2.
\end{equation}
By combining \eqref{pr4.5 1} and \eqref{core 4} with the orthogonality condition \eqref{pro2},  the desired inequality \eqref{core 1} follows.    \end{proof} 
  We next establish the main regularity result.
\begin{proposition}\label{regularity}
     For the $\uu^\pm$ obtained in Corollary~\ref{weak con}, we have $$\uu^\pm\in L^{\infty}\big(0,T;H^1(\Omega_t^\pm;N^\pm)\big)$$ 
    and there exists a constant $C=C(\g)>0$ independent of $\e
    $ such that
    \begin{equation}\label{claim 4.5.2}
       \esssup_{t\in(0,T)} \sum_{\pm}\int_{\o^\pm_t}|\nabla\uu^\pm(\cdot,t)|^2dx\leq C.
    \end{equation}
\end{proposition}
\begin{proof}   
   Let $\e_k\downarrow0$ be the subsequence along which the  conclusions regarding $\uek$ and $\uu^\pm$ in  Corollary~\ref{weak con} hold. For each $k$, let $\delta=\delta_k=\e_k^{\frac{1}{6}}$ as in Proposition~\ref{level set estimate}, so that $a^\pm_{\e_k,\delta_k}\leq C\delta_k^{-2}\e_k^{\frac{1}{2}}=C\e_k^\frac{1}{6}\downarrow0$.  We define the sets 
    \begin{align}
        \o_t^{k,-}&=\{\psi_{\e_k}< a^-_{\e_k,\delta_k}(t)\},\\
         \o_t^{k,+}&=\{\psi_{\e_k}>c_F-a^+_{\e_k,\delta_k}(t)\}.
    \end{align}
      Since $a^\pm_{\e_k,\delta_k}\downarrow0$, there exists $K_1\in\mathbb{N}^+$  independent of $t$ such that for all $k\geq K_1$,
      \begin{equation}\label{well-defined}
          \begin{aligned}
              \psi_{\e_k}> c_F-a^+_{\e_k,\delta_k}\quad&\Rightarrow\quad \uek\in B_{\delta_1}(N^+),\\
          \psi_{\e_k}<a^-_{\e_k,\delta_k}\quad&\Rightarrow \quad\uek\in B_{\delta_1}(N^-),
          \end{aligned}
      \end{equation}      
      with $\delta_1:=\min\{\frac{1}{2C(N)},\delta_N\}$, where $C(N)$ is from Proposition~\ref{pr core} and $\delta_N$ is defined in~\eqref{delta_0}. 
      
      Define the  projection truncation function (to project onto $N$)
    \begin{equation}
        \vv_k(x,t)=\sum_\pm P_N\big(\uu_{\e_k}(x,t)\big)\chi_{\o_t^{k,\pm}}\quad\text{ for } k\geq K_1,
    \end{equation}
    which is well-defined due to~\eqref{well-defined}.
    
    By Proposition~\ref{level set estimate} and the regularity $\uek\in C\big([0,T],H^1(\o)\big)$ (cf. Lemma \ref{ue regularity}), we have 
\begin{align}
\vv_k
    \in L^\infty\big([0,T]&,SBV(\o)\big),\\
    \nabla^a\vv_k=\sum_{\pm}\nabla \big(P_N(\uek)\big)\chi_{\otknp},&\quad
    J_{\vv_k}\overset{\hn}=\bigcup_\pm\p^*\otknp\cap\o.\label{4.51}
\end{align}  
      Using Proposition~\ref{pr core} with $\delta=\delta_1\leq\frac{1}{2C(N)}$ and the estimate~\eqref{cor3.4 0.3}, we obtain a constant $C > 0$ such that
      \begin{equation}\label{tange bound}
           \sup_{t\in[0,T]}\int_{\otknp}\Big|\nabla \big(P_N(\uek)\big)\Big|^2dx\leq C\quad\forall k\geq K_1.
      \end{equation}     
 Hence, combining~\eqref{4.51} and~\eqref{tange bound} with Proposition~\ref{level set estimate}, we conclude that
 \begin{equation}
\sup_{k\geq K_1}\esssup_{t\in(0,T)}\Big(\|\vv_k\|_{L^\infty(\o)}+\|\nabla^a\vv_k\|_{L^2(\o)}+\hn(J_{\vv_k})\Big)<\infty.\label{bounded}
 \end{equation}
  By the compactness of $BV$ spaces with uniformly bounded total variation (cf. \cite[Theorem 3.23]{ambrosio2000functions}) and Lemma~\ref{sbv lemma},  there exists a function $\vv\in L^\infty\big(0,T;BV(\o)\big)$ such that, up to a subsequence,
    \begin{equation}        
\vv_k\xrightharpoonup{k\to\infty}\vv\in SBV(\o)\quad\text{ weakly-star in } L^\infty\big(0,T;BV(\o)\big).\label{v weakly star}
    \end{equation}
    
   We proceed to identify the limit $\vv$. The weak-star convergence~\eqref{v weakly star}, the estimate~\eqref{bulk error} and the pointwise convergence~\eqref{limit 4} imply the following facts:  for a.e. $t\in(0,T)$, up to a subsequence,
    \begin{equation}\label{pointwise}
    \begin{aligned}
       \vv_k(\cdot,t)&\xrightarrow{k\to\infty}\vv(\cdot,t) \quad\text{ a.e. in } \o,\\       \chi_{\o_t^{k,\pm}}&\xrightarrow{k\to\infty}\chi_{\o_t^{\pm}}\quad\text{ a.e. in } \o,\\
       \uek(\cdot,t)&\xrightarrow{k\to\infty}\uu^\pm(\cdot,t) \quad\text{ a.e. in } \o_t^\pm.
    \end{aligned}              
    \end{equation}
    Fix such $t$, and let $K\subset\o_t^\pm$ be a compact subset. By~\eqref{well-defined}, for all $k \geq K_1$ and $x \in \Omega_t^{k,\pm}$, we have $\uu_{\e_k}(x,t) \in B_{\delta_N}(N^\pm)$, therefore, $F(\uek)=\dd_N^2(\uek)=|P_N\uek-\uek|^2$ (recall the definition of $F$ from~\eqref{f}).  Applying Fatou’s lemma together with the Cauchy–Schwarz inequality,~\eqref{pointwise} gives
    \begin{equation}
        \begin{aligned}
            \int_K\big|\vv-\uu^\pm\big|dx\leq&\liminf_{k\to\infty}\int_K|P_N\uek-\uek|\chi_{\otknp}dx\\
                          \leq&\liminf_{k\to\infty}|K|^\frac{1}{2}\bigg(\int_KF(\uek)dx\bigg)^{\frac{1}{2}}\\
                          \overset{\eqref{estimate b}}=&0,
        \end{aligned}
    \end{equation}
 yielding $\vv= \sum_{\pm}\uu^\pm\chi_{\o_t^{\pm}}$ a.e. in $\o\times(0,T)$. 

Since $\nabla\uu^\pm\in L^\infty\big(0,T;L^2_{\mathrm{loc}}(\o^\pm_t)\big)$, we deduce $\hn\big(J_\vv\cap\o^\pm_t\big)=0$ for a.e. $t\in(0,T)$. Hence, by the structure theorem of $SBV$ functions (cf. \cite[Section 4.1]{ambrosio2000functions}),  $\nabla^a \vv = \nabla \uu^\pm$ on $\Omega_t^\pm$.  Finally, by~\eqref{inequality 11} in Lemma~\ref{sbv lemma}, we obtain
    \begin{equation}
\sum_{\pm}\int_{\o^\pm_t}|\nabla\uu^\pm(\cdot,t)|^2dx\leq\liminf_{k\to\infty}\int_{\o}|\nabla^a\vv_k(\cdot,t)|^2dx.
    \end{equation}
This combined with the uniform bound~\eqref{bounded}, implies the desired regularity~\eqref{claim 4.5.2}, thus completing the proof. 

\end{proof}

\subsection{Minimal pair conditions}
Once the regularity of the limit functions $\uu^\pm$ is established, we are in a position to adopt the strategy developed by Liu \cite[Section 5]{liu2022phase} with some adaptations for the case of boundary contact.

We begin by introducing the semi-distance,
\begin{equation}
     \dd_F^*(\uu_+,\uu_-):=\inf_{\substack{\xi(\pm1)=\uu_\pm\\\xi\in H^1\big((-1,1),\rk\big)}}\int_{-1}^{1}|\xi^\prime(t)|\sqrt{2F\b(\xi(t)\b)}dt,
\end{equation}
which induces a semi-distance between closed sets $S_\pm\subset\rk$ by
   \begin{equation}
       \dd_F^*(S_+,S_-):=\inf_{\uu_\pm\in S_\pm}\dd_F^*(\uu_+,\uu_-).
   \end{equation}
For any $S_\pm\subset N^\pm$, we define   
\begin{equation}
    \mathcal{N}(S_\pm,\delta):=\bigcup_{\uu_\pm\in S_\pm}\mathcal{N}(\uu_\pm,\delta),
\end{equation}
where $\mathcal{N}(\uu_\pm,\delta)$ is the normal sphere induced by the metric $\dd_F^*$ given by
\begin{equation}
   \mathcal{N}(\uu_\pm,\delta):=\{\uu\in\rk: \dd_F^*(\uu,\uu_\pm)=\delta,P_N(\uu)=\uu_\pm\}.
\end{equation}
We then define a function $\kappa:N^+\times N^-\times[0,1]\rightarrow \mathbb{R}$ as
\begin{equation}
\kappa(\uu_+,\uu_-,\rho):=\dd_F^*\bigg(\mathcal{N}\Big(\overline{B_\rho(\uu_+)}\cap N^+,\delta_N\Big) ,\mathcal{N}\Big(\overline{B_\rho(\uu_-)}\cap N^-,\delta_N\Big)\Bigg)+2{\delta_N}-c_F,
\end{equation}
where $\delta_N$ is given in~\eqref{delta_0} and $c_F$ is defined in~\eqref{cf}. 

The following properties are recorded in Liu \cite{liu2022phase}:
\begin{lemma}\cite[Lemma 5.1]{liu2022phase}\label{lemma 4.7}
     The function $\kappa(\uu_+,\uu_-,\rho)$ is non-negative. Moreover, $\kappa(\uu_+,\uu_-,\rho)=0$ if and only if $(\uu_+,\uu_-)\in N^+\times N^-$ is a minimal pair.
\end{lemma}
\begin{lemma}\cite[Lemma 5.2]{liu2022phase}\label{lemma 4.8}
    There exist constants $C_0=C_0(N)$ and $\rho_0=\rho_0(N)$ which depend only on the geometry of $N$ such the following holds: if
    \begin{equation}
        \gamma\in H^1([-\delta,\delta],\rk)\text{ with }\gamma(\pm\delta)\in B_{{\delta_N}} (N^\pm), \text{ and } \rho\in(0,\rho_0),
    \end{equation}
    then
    \begin{equation}
        \int_{-\delta}^{\delta}\Big(\frac{1}{2}|\gamma^\prime|^2+\frac{1}{\e^2}F(\gamma)-\frac{1}{\e}\big(\dd_F(\gamma)\big)^\prime\Big)\geq\frac{\min\Big\{C_0\rho^2,\kappa\Big(P_N\big(\gamma(\delta)\big),P_N\big(\gamma(-\delta)\big),\rho\Big)\Big\}}{\max\{\e,\delta\}}.
    \end{equation}
\end{lemma}
We proceed via a contradiction argument. Let $\e_k\downarrow0$ be the subsequence along which the conclusions regarding $\uek$ and $\uu^\pm$ in Corollary~\ref{weak con} and~\ref{regularity} hold. Assume that there exist $t\in[0,T]$, a constant $\alpha>0$, and a compact subset $E_t^*\subset\Gamma^*_t:=\Gamma_t\cap\{x\in\o:\dist(x,\p\o)>\alpha\}$ such that 
\begin{align*}
    \uu^\pm(\cdot,t)\in H^1(\o_t^\pm),\quad&\hn(E_t^*)\geq\alpha,\\
    \forall p\in E_t^*,\quad\big(\uu^+(p),\uu^-(p)\big)\subset& N^+\times N^- \text{ is not a minimal pair. }
\end{align*}
The following argument is for the fixed $t$.

Let $\nn := \nn_{\Gamma_t}$ denote the unit normal to $\Gamma_t$.
By the definition of $\g_t^*$, as in~\eqref{diffeo}, there exists a constant $\delta_2>0$ such that the mapping
\begin{equation}
    \Psi_{\Gamma_t^*}:\Gamma^*_t\times(-\delta_2,\delta_2)\rightarrow\o,\quad(x,s)\mapsto x+s\nn(x)
\end{equation}
is a $C^2$ diffeomorphism  with bounds $|\nabla\Psi_{{\Gamma_t^*}}|,|\nabla\Psi_{{\Gamma_t^*}}^{-1}|\leq
C$. Using the strong convergence~\eqref{limit 3} and the trace property of Sobolev functions (cf. \cite[Theorem 5.6]{evans2018measure}), there exists a null set $N\subset[0,{\delta_2}]$ such that
\begin{align}    &\uu_{\e_k}\big(p\pm\delta\nn(p)\big)\xrightarrow{k\to\infty}\uu^\pm\big(p\pm\delta\nn(p)\big) &&\text{ strongly in } L^2(E_t^*
    ) \quad \forall\delta\notin N,\\
&\uu^\pm\big(p\pm\delta\nn(p)\big)\xrightarrow{\delta\downarrow0,\delta\notin N}\uu^\pm(p) &&\text{ strongly in } L^2(E_t^*).
\end{align}
A diagonal argument yields a sequence $\delta_k\downarrow0$ such that
\begin{equation}    \uu_{\e_k}\big(p\pm\delta_k\nn(p)\big)\xrightarrow{k\to\infty}\uu^\pm(p) \quad\text{ strongly in } L^2(E_t^*
    ) .
\end{equation}
Next, we consider the orthogonal decomposition of the gradient  $\nabla=\nabla_\Gamma+\nn\p_\nn$  in $\Psi_{\Gamma_t^*}\big(\g_t^*\times(-\delta_2,\delta_2)\big)$.
Then, from the estimate~\eqref{E estimate}, we have
\begin{align}\label{mini a}
    \int_{\Psi_{\Gamma_t^*}\b(\g_t^*\times(-\delta_2,\delta_2)\b)}\bigg(\frac{\e
    }{2}|\p_\nn\uek|^2+\frac{1}{\e}F(\uek)-\bxi\cdot\nabla_\g\pek-\bxi\cdot\nn\p_\nn\pek \bigg)dx\leq C\e.
\end{align}
From the regularity of $\bxi$~\eqref{eq_regularityXi} and the consistency condition~\eqref{eq_consistencyProperty}, a Taylor expansion yields for $x\in\Psi_{\g_t^*}\big(\g_t^*\times(-\delta_2,\delta_2)\big)$,
\begin{equation}\label{pair 1}
\begin{aligned}
     |\bxi\cdot\nabla_\g\pek|\leq C\dd_\g^2|\nabla\pe|,\quad 
      |(\nn-\bxi)\cdot\nn\p_\nn\pek|\leq C\dd_\g^2|\nabla\pe|.
\end{aligned}   
\end{equation}
Applying this together with ~\eqref{E estimate} and the coercivity estimate~\eqref{energy f}, we obtain
\begin{equation}
      \int_{\Psi_{\Gamma_t^*}\big(\g_t^*\times(-\delta_2,\delta_2)\big)}\Big( \big|\bxi\cdot\nabla_\g\pek\big| dx+|\p_\nn\pek-\bxi\cdot\nn\p_n\pek|\Big)dx\leq C\e.\label{mini 1}
\end{equation}
Combining~\eqref{mini 1} with~\eqref{mini a}, using the change of variables $x = p + s \nn(p)$ along with the area formula and the bounds $|\nabla\Psi_{\g_t^*}|,|\nabla\Psi_{\g_t^*}^{-1}|\leq
C$, we arrive at
\begin{equation}
    \int_{E_t^*}\underbrace{\int_{-\delta_k}^{\delta_k}\lef\frac{1}{2}|\p_s\uu_{\e_k}|^2+\frac{1}{\e^2}F(\uu_{\e_k})-\frac{1}{\e}\p_s\big(\dd_F(\uu_{\e_k})\big)\rig d\tau}_{:=X_k(p)}dS(p)< \infty.\label{contradition 1}
\end{equation}

On the other hand,  since  $\big(\uu^+(p), \uu^-(p)\big)$ is not a minimal pair for $p \in E_t^*$, Lemma~\ref{lemma 4.7} gives
\begin{equation}\label{4.74}
    \kappa\big(\uu^+(p),\uu^-(p),0\big)>0\quad\forall p\in E_t^*.
\end{equation}
By Egorov's theorem, there exists a compact subset $E_t\subset E_t^*$ with $\hn(E_t)>\frac{\alpha}{2}$ such that
\begin{equation}
    \uu_{\e_k}\big(p+\delta_k\nn(p)\big)\xrightarrow{k\to\infty}\uu^\pm(p) \text{ uniformly on } E_t.
\end{equation}
This implies that there exists $K\in \mathbb{N}^+$, such that
\begin{equation}
    \uek\big(p\pm\delta_k\nn(p)\big)\in B_{\delta_N}(N^\pm)\quad \forall k\geq K,p\in E_t.
\end{equation}
Moreover, since $\uu_{\e_k}\big(p\pm\delta_kn(p)\big)$ are smooth and converge uniformly on $E_t\subset\Gamma^*_t$, the functions $\uu^\pm$ are continuous on $E_t$. Thus, the mapping
\begin{equation}
  \kappa\big(\uu^+(\cdot),\uu^-(\cdot),\cdot\big):E_t\times[0,1]\rightarrow[0,\infty)
\end{equation}
is continuous, and  from \eqref{4.74}, there exist $\rho\in(0,\rho_0)$ and $\beta>0$ such that
\begin{equation}
    \inf_{p\in E_t}\kappa\big(\uu^+(p),\uu^-(p),\rho\big)=2\beta.
\end{equation}
By the uniform convergence of $\uek$, possibly increasing $K$ if necessary, we have
\begin{equation}
     \inf_{p\in E_t}\kappa\bigg(P_N\Big(\uu_{\e_k}\big(p+\delta_k\nn(p)\big)\Big),P_N\Big(\uu_{\e_k}\big(p-\delta_k\nn(p)\big)\Big),\rho\bigg)\geq\beta\quad\forall k\geq K.
\end{equation}
Finally, by Lemma~\ref{lemma 4.8}, we deduce that
\begin{equation}
    X_k(p)\geq\frac{\min\{C_0\rho^2,\beta\}}{\max\{\e_k,\delta_k\}}\quad\forall p\in E_t, k\geq K, 
\end{equation}
contradicting the finiteness of the integral ~\eqref{contradition 1}.

This contradiction completes the argument.
\subsection{Harmonic heat flows in the bulk}
In this subsection, we prove that $\uu^\pm$ solves the harmonic heat flow into $N^\pm$ as described in Theorem~\ref{theorem 2}. Although this result is standard, we include the details here for completeness. We refer to \cite[Section~4]{liu2022phase}, \cite[Section~7.6]{lin2008analysis}, and \cite{chen1989existence} for background. In the following, we argue for $\uu^+$; this argument for $\uu^-$ is symmetric.

Let $\mathcal{P}_R(z_0):=B_R(x_0)\times(t_0-R^2,t_0)$ denote  the parabolic cylinder of size $R$ centered at $z_0=(x_0,t_0)$. Fix $\delta>0$ and a cylinder $$\mathcal{P}_{4R_0}(z_0)\subset\bigcup_{t\in(0,T)}\b(\otp\setminus B_\delta(\Gamma_t)\b)\times\{t\}.$$ Let $\{\ue\}$ be the weak solution from Theorem~\ref{theorem 1} and recall that $\|\ue\|_{L^\infty}\leq R_N$ from Lemma~\ref{ue regularity}. Using the estimate~\eqref{estimate b}, we define
\begin{equation}
    E_\delta:=\sup_{t\in[0,T]}\int_{\o_t^+\setminus B_\delta(\Gamma_t)}\bigg(\frac{|\nabla\ue|^2}{2}+\frac{1}{\e}F(\ue)\bigg)dx\leq C\delta^{-2}.
\end{equation}
 
\begin{lemma}\cite[Lemma 4.3]{chen1989existence}\label{bocher}
    There exists $c=c(N)>0$ independent of $\e$ such that for
    \begin{equation}
        e_\e(\ue):=\frac{|\nabla\ue|^2}{2}+\frac{1}{\e^2}F(\ue)
    \end{equation}
    satisfies  
    \begin{equation}
        (\p_t-\Delta)e_\e(\ue)+\e^{-4}\dd_N^2(\ue){f^{\prime}}^2(\ue)\leq c\big(1+e_\e(\ue)\big)e_\e(\ue) \text{ in } \mathcal{P}_{4R_0}(z_0).
    \end{equation}
\end{lemma}
Fix any $z_1=(x_1,t_1)\in \overline{\mathcal{P}_{2R_0}(z_0)}$ and consider the back heat kernel centered at this point,
\begin{equation}
    G_{z_1}(x,t):=\frac{1}{(4\pi|t-t_1|)^\frac{d}{2}}e^{-\frac{|x-x_1|^2}{4|t-t_1|}}\quad\text{ for }x\in\R^n,t< t_1.
\end{equation}
 Let $\eta_0\in C^\infty_c\big(B_{2R_0}(0)\big)$ be a cut-off function such that ${\eta_0}\equiv1$ in $B_{R_0}(0)$. Define $\eta_1(x)=\eta_0(x-x_1)\in C^\infty_c\big(B_{2R_0}(x_1)\big)$. For any $r<R_0$, set
\begin{equation}
    \begin{aligned}
    \Phi_r(\ue, z_1)&: = r^2 \int_{\mathbb{R}^n \times \{t = t_1 - r^2\}} \eta_1^2 G_{z_1} e_\e(\ue) dx,\\
      \Psi_r(\mathbf{u}_\e, z_1)&: = \int_{T_r(z_1)} \eta_1^2 G_{z_1} e_\e(\ue)dxdt,
    \end{aligned}
\end{equation}
where $T_r(z_1)=\R^n\times(t_1-4r^2,t_1-r^2)$.
\begin{lemma}\cite[Lemma 4.2]{chen1989existence}
    For $0<r_1\leq r_2<R_0$, we have
    \begin{equation}
        \begin{aligned}
            \Psi_{r_1}(\ue, z_1) \leq \Psi_{r_2}(\ue, z_1) + (r_2 - r_1) C E_\delta,\\
            \Phi_{r_1}(\ue, z_1) \leq \Phi_{r_2}(\ue, z_1) + (r_2 - r_1) C E_\delta.
        \end{aligned}
    \end{equation}
\end{lemma}
\begin{lemma}\cite[Lemma 4.4]{chen1989existence}\label{small energy}
    There exist $\overline{\e},\overline{\delta},c>0$ depending only on $N$ (independent of $z_1$) such that for any $z_1=(x_1,t_1)\in \overline{\mathcal{P}_{2R_0}(z_0)}$ if 
    \begin{equation}
          \Psi_{r}(\ue, z_1)\leq\overline{\e}^2\quad\text{ for some } r\in(0,R_0),
    \end{equation}
    then 
    \begin{equation}
        \sup_{z\in \mathcal{P}_{\overline{\delta}r}(z_1)}e_\e(\ue)(z)\leq c(\overline{\delta}r)^{-2}.
    \end{equation}
\end{lemma}

Let $\e_k\downarrow0$ be the subsequence along which the conclusions regarding $\uek$ and $\uu^\pm$ in Corollary~\ref{weak con} and~\ref{regularity} hold. Noticing the arbitrariness of the choice of $z_1\in\overline{\mathcal{P}_{2R_0}(z_0)}$, we define the energy concentration set
\begin{equation}
    \mathcal{S}:=\bigcap_{0<r<R_0}\bigg\{z\in\overline{\mathcal{P}_{2R_0}(z_0)}:\liminf_{k\to\infty}  \Psi_r(\uek, z)\geq\overline{\e}^2\bigg\},
\end{equation}
where $\overline{\e}$ is given in Lemma~\ref{small energy}.
It can be shown that $\mathcal{S}$ is closed and has locally finite $n$-dimensional parabolic Hausdorff measure, with $\mathcal{P}^n(\mathcal{S})\leq CE_\delta$ (see \cite[Section 7.6]{lin2008analysis} for the proof). Fix any $(r,z^\prime)\in(0,R_0)\times(\mathcal{P}_{2R_0}(z_0)\setminus\mathcal{S})$, by the construction of $\mathcal{S}$, there exists $K=K(r,z^\prime)\in\mathbb{N}^+$ such that  
$$ \Psi_r(\uek, z^\prime)<\overline{\e}^2\quad\forall k\geq K.$$ 
Then, by Lemma~\ref{small energy}, for such $(r,z^\prime)$, we have
\begin{equation}
      \sup_{z\in \mathcal{P}_{\overline{\delta}r}(z^\prime)}e_{\e_k}(\uek)(z)\leq c(\overline{\delta}r)^{-2} \quad\forall k\geq K.
\end{equation}
For any $\{\overline{r}\}\times\mathcal{P}_{2\overline{r}}(\overline{z})\subset (0,R_0)\times(\mathcal{P}_{2R_0}(z_0)\setminus\mathcal{S})$, by applying the above argument and a covering argument with the family of sets $\{\mathcal{P}_{\overline{r}}(z):\mathcal{P}_{\overline{r}}(z)\subset \mathcal{P}_{2R_0}(z_0)\setminus\mathcal{S}\}$, we conclude that there exists $K=K(\overline{r},\overline{z})\in\mathbb{N}^+$ such that 
\begin{equation}\label{small energy regularity}
    \sup_{z\in \mathcal{P}_{2\overline{\delta}\overline{r}}(\overline{z})}\Big(|\nabla\uek|^2 +{\e}_k^{-2}F\big(\uek(z)\big)\Big)\leq c(\overline{\delta}\overline{r})^{-2}\quad\forall k\geq K.
\end{equation}
This uniform bound further implies that $\uek$ is uniformly bounded in $C^\frac{1}{2}\big(\mathcal{P}_{\overline{\delta}\overline{r}}(\overline{z})\big)$ (see \cite[Section 7.6]{lin2008analysis} for details). After passing to a subsequence, we obtain
\begin{align}
    \uek&\xrightarrow{k\to\infty}\uu^+\quad\text{ strongly in }C^0\big(\mathcal{P}_{\overline{\delta}\overline{r}}(\overline{z})\big), \label{strong 1}\\\nabla\uek&\xrightarrow{k\to\infty}\nabla\uu^+\quad\text{ weakly-star in }L^\infty\big(\mathcal{P}_{\overline{\delta}\overline{r}}(\overline{z})\big).
\end{align}
From the bulk estimate~\eqref{unweighted}, by increasing $K$ further if necessary, we may assume that, recalling the definition of $f$ from~\eqref{f},
\begin{equation}
    f^\prime\big(\dd_N^2(\uek)\big)=1\quad\forall k\geq K.
\end{equation}
Combining Lemma~\ref{bocher}  with  the bound~\eqref{small energy regularity}, we get
\begin{equation}
    (\p_t-\Delta)e_{\e_k}(\uek)+\e_k^{-4}\dd^2_N(\uek)\leq C \quad\forall k\geq K.\label{ineq 1}
\end{equation}
Choose a nonnegative test function $\phi\in C^\infty_c\big(\mathcal{P}_{\overline{\delta}\overline{r}}(\overline{z})\big)$. Testing~\eqref{ineq 1} against $\phi$ and  integrating by parts, we obtain
\begin{equation}
    \int_{\mathcal{P}_{\overline{\delta}\overline{r}}(\overline{z})}\e_k^{-4}\dd_N^2(\uek)\phi dxdt
    \leq \int_{\mathcal{P}_{\overline{\delta}\overline{r}}(\overline{z})}\big(|\p_t\phi+\Delta\phi|e_{\e_k}(\uek)+C\phi\big)dxdt\leq C(\phi).
\end{equation}
Moreover, $\uek$ satisfies the equation
\begin{equation}
    \p_t\uek=\Delta\uek-\e_k^{-2}D(\dd_N^2)(\uek).\label{equation 1}
\end{equation}
It follows that $\e_k^{-2}\dn(\uek)$ is uniformly bounded in $L^2_{\mathrm{loc}}\big(\mathcal{P}_{\overline{\delta}\overline{r}}(\overline{z})\big)$. Hence, $(\p_t-\Delta)\uek$ is uniformly bounded in $L^2_{\mathrm{loc}}\big(\mathcal{P}_{\overline{\delta}\overline{r}}(\overline{z})\big)$, and $\p_t\uek,\nabla^2\uek$ are also uniformly bounded in $L^2_{\mathrm{loc}}\big(\mathcal{P}_{\overline{\delta}\overline{r}}(\overline{z})\big)$. By the  weak compactness, after passing to a further subsequence,
\begin{align}
      (\p_t\uek,\nabla^2\uek)&\xrightarrow{k\to\infty}(\p_t\uu^+,\nabla^2\uu^+) \quad \text{ weakly in }L^2_{\mathrm{loc}}\big(\mathcal{P}_{\overline{\delta}\overline{r}}(\overline{z})\big).\label{weak 1}
\end{align}
There exists a unit vector field $\bnu_N^k\perp T_{P_N\uek}N$ such that
\begin{equation}
    D(\dd_N^2)(\uek)=2\dn(\uek)\bnu_{N}^k.
\end{equation}
For any vector field $\psi\in L^2_{\mathrm{loc}}\big(\mathcal{P}_{\overline{\delta}r}(z) \big)$ with $\psi(z)\in T_{P_N\uu(z)}N$, from~\eqref{strong 1}, we have
\begin{align}
    \lim_{k\to\infty}\int_{\mathcal{P}_{\overline{\delta}\overline{r}}(\overline{z})}\e_k^{-2}D(\dd_N^2)(\uek)\cdot \psi=0.
\end{align}
Combining this with~\eqref{equation 1} and~\eqref{weak 1},  we conclude that $(\p_t-\Delta)\uu^+\perp T_{\uu^+}N$ a.e. in  $\mathcal{P}_{\overline{\delta}\overline{r}}(\overline{z})$.

Finally, by standard parabolic regularity theory, $\uu^+$ is smooth in $\mathcal{P}_{\overline{\delta}\overline{r}}(\overline{z})$. Thus, $\uu^+$ is a classical solution away from the singular set $\mathcal{S}$. By a capacity argument, noticing the arbitrariness of the choice of  $\mathcal{P}_{2\overline{\delta}\overline{r}}(\overline{z})$ , we conclude that $\uu^+$ extends as a weak solution in $\mathcal{P}_{2R_0}(z_0)$.

\section*{Acknowledgments}
The author gratefully acknowledges Professor Fanghua Lin  and Professor  Yuning Liu for insightful discussions, and expresses deep appreciation to supervisor Professor Yaguang Wang for his guidance. 
\printbibliography

\appendix
\section{An outline of the proof of Lemma~\ref{ue regularity}}\label{proof lem 2.6}
Let $T>0$ be fixed, $N\in\mathbb{N}^+$ and $\tau=\tau(N):=\frac{T}{N}$. Define $\uu_N^0:=\uu_{\e,0}\in H^1(\o)$, satisfying the setting in Lemma \ref{ue regularity}. For $k=1,\cdots,N$, we construct $\uu_N^k$ inductively : given $\uu_N^{k-1}$, then let $\uu_N^k$ be a minimizer of the functional 
\begin{equation}
    E_k(\cdot):H^1(\o)\rightarrow[0,\infty]:\uu\rightarrow E_\e(\uu)+\frac{1}{2\tau}\|\uu-\uu_N^{k-1}\|^2_{L^2(\o)}.
\end{equation}
The existence of a minimizer follows from  the  growth condition \eqref{convexity} on $F$.
The associated Euler-Lagrange equation is given by (see \cite[Lemma 3.5]{garcke2003cahn}): for all test functions $\xi\in H^1(\o)\cap L^\infty(\o)$,
\begin{equation}\label{euler-lagr}
\begin{aligned}
    \int_\o\bigg(\e\nabla\uu_N^k:\nabla\xi+\frac{\uu_N^k-\uu_N^{k-1}}{\tau}\cdot\xi+\frac{1}{\e}&DF(\uu_N^k)\cdot\xi\bigg)dx\\&+\int_{\p\o}\nabla\sigma(\uu_N^k)\cdot\xi d\hn=0.
\end{aligned}  
\end{equation}
We claim that the minimizers \(\uu_N^k\) can be chosen such that \(\|\uu_N^k\|_{L^\infty(\Omega)}\leq R_N\). To establish this claim, we rely on the following lemma.
\begin{lemma}\label{radial trun}
     Fix any constant $R>0$, define the radial truncation $ T_{R}:\rk\rightarrow\rk$ by
    \begin{equation}
        T_{R}(\uu)=\begin{cases}
            \uu&\text{ if }|\uu|\leq R,\\
            R\cdot\frac{\uu}{|\uu|}&\text{ if } |\uu|>R.
        \end{cases}
    \end{equation}
   Then, for all \(\uu,\vv\in\mathbb{R}^k\), the following properties hold:
    \begin{subequations}
        \begin{align}            
            |T_{R}(\uu)|&\leq R, \label{a3}\\      
             |T_{R}(\uu)-T_{R}(\vv)|&\leq|\uu-\vv|,\label{a4}
        \end{align}
    \end{subequations}
    and for all \(\uu\in H^1(\o)\), 
    \begin{equation}\label{a5}
        \|\nabla\big(T_{R}(\uu)\big)\|_{L^2(\o)}\leq\|\nabla\uu\|_{L^2(\o)}.
    \end{equation}
\end{lemma}
   
\begin{proof}
    The inequality~\eqref{a3} follows immediately from the definition of \(T_{R}\). The inequality~\eqref{a5} can be derived as a corollary of~\eqref{a4} by applying the generalized chain rule of Ambrosio and Dal Maso. We prove the inequality~\eqref{a4} via a case-by-case analysis:
  
  Case 1: $|\uu|,|\vv|\leq R$. In this case, \(T_{R}(\uu)=\uu\) and \(T_{R}(\vv)=\vv\), so the inequality is trivially satisfied.

 Case 2: $|\uu|> R_2$ and $|\vv|>R$. We have
    \begin{equation}
        \begin{aligned}
            |T_{R}(\uu)-T_{R}(\vv)|=&R\bigg|\frac{\uu}{|\uu|}-\frac{\vv}{|\vv|}\bigg|\\
            \leq& \frac{R_2}{\min\{|\uu|,|\vv|\}}|\uu-\vv|\\
            \leq&|\uu-\vv|.            
        \end{aligned}
    \end{equation}

 Case 3: $|\uu|\leq R$ and $|\vv|>R$. Consider the triangle formed by the points \(\uu,\vv,T_{R_2}(\vv)\) in \(\mathbb{R}^k\). By geometric considerations, the angle between the vectors \(\uu - T_{R_2}(\vv)\) and \(\vv - T_{R_2}(\vv)\) is at least \(90^{\circ}\). This implies that
    \begin{equation}       
           |T_{R}(\uu)-T_{R}(\vv)|=|\uu-T_{R}(\vv)|\leq|\uu-\vv|.
    \end{equation}

\end{proof}
By this Lemma and assumptions~\eqref{F 2} and~\eqref{sigma 1} on  $F$ and \(\sigma\), for any \(\uu\in H^1(\Omega)\), we have
\begin{equation}
    E_k\b(T_{R_N}(\uu)\b)\leq E_k(\uu)\quad \text{ if } \|\uu_N^{k-1}\|_{L^\infty(\o)}\leq R_N.
\end{equation}
Next, we define the piecewise-linear extension $\uu_N$ of the sequence $\{\uu_N^k\}$. Specifically,  for $t=\lambda(k-1)\tau+(1-\lambda)k\tau$, where $\lambda\in[0,1]$ and $k=0,\cdots,N$, we set
\begin{equation}
    \uu_N(t):=\lambda\uu_{N}^{k-1}+(1-\lambda)\uu_N^k.
\end{equation}
Following \cite[Appendix A]{hensel2022convergence} and using  assumptions on $F$ and \(\sigma\), we can show that there exists a subsequence of \(\{\uu_N\}\) (not relabel) such that
\begin{equation}
    \uu_N\rightarrow\uu\text{ in } L^2\b(0,T,H^1(\o)\b)\cap H^1\b(0,T,L^2(\o)\b).
\end{equation}
 We prove these statements in Lemma~\ref{ue regularity} (also cf. \cite[Appendix A]{hensel2022convergence}).
\begin{itemize}
    \item[(1)]  \textit{Uniform boundedness}. It follows directly from the choice of $\uu_N^k$.
    \item[(2)]  \textit{Variational formulation and uniqueness}. By leveraging the above convergence properties, we can take the limit in the associated Euler-Lagrange equation. The uniqueness of weak solutions can be established via a Gronwall-type argument combined with the $L^\infty$-bound.
    \item[(3)] \textit{High regularity and boundary conditions}.
             Since \(\|\uu\|_{L^\infty(\Omega)}\leq R_N\), and under the assumptions on $F$ and \(\sigma\), we can improve the regularity of \(\uu\) following the same procedure as in \cite[Appendix A]{hensel2022convergence}.
    \item[(4)] \textit{Interior regularity}.  Due to the specific form of the function $F$, we can apply the standard bootstrapping method and interior regularity estimates for parabolic equations to obtain the desired interior smoothness.
    
    \item[(5)]  \textit{Energy dissipation identity}. It follows from the high regularity by a standard mollification argument.
\end{itemize}

\section{Construction of well-prepared initial data}\label{constru}
\vspace{0.5em}
\noindent\textbf{Initial data construction.}  
Recall the function $f$ defined in~\eqref{f}. We introduce the associated scalar potential
\begin{equation}
\widetilde{F}(\lambda):= \begin{cases}f\left(\left(\frac{\dist_N}{2}+\lambda\right)^{2}\right) & \text { if } \lambda \leq 0, \\
f\left(\left(\frac{\dist_N}{2}-\lambda\right)^{2}\right)& \text { if } \lambda \geq 0.\end{cases}\label{1.2.6}
\end{equation} 
Consider the unique solution to the following boundary value problem
\begin{equation}\label{1.2.10-1}
 \begin{cases}
\alpha^{\prime }(t)=\sqrt{2\widetilde{F}\b(\alpha(t)\b)} \quad t \in \mathbb{R},\\[2mm]
\alpha(\pm \infty)= \pm\frac{\dist_N}{2}.
\end{cases}
\end{equation}
The function \(\alpha(t)\) (understood as the optimal profile in phase transition) satisfies the following properties (cf. \cite{fonseca1989gradient,sternberg1988effect}):
\begin{itemize}
    \item[(1)] It is odd and strictly monotone increasing, and 
    \begin{equation} 
 \alpha \in C^{\infty}\left(\mathbb{R}, \left(-\frac{\dist_N}{2},\frac{\dist_N}{2}\right)\right) \text{ with } \alpha(0)=0.\label{1.2.9}
\end{equation}
    \item[(2)] There exist constants $C_1,C_2>0$ such that 
    \begin{equation}\label{1.2.10}
 \begin{cases} 
\left|\alpha^{\prime}(t)\right|+\left|\alpha(t)+\frac{\dist_N}{2}\right| \leq C_{1} e^{C_{2} t}  &\text { when } t <0, \\
\left|\alpha^{\prime}(t)\right|+\left|\alpha(t)-\frac{\dist_N}{2}\right| \leq C_{1} e^{-C_{2} t}  &\text { when } t>0.
\end{cases}
\end{equation}

\end{itemize}


\begin{lemma}
Define the truncated function $\alpha_{2L}$ of $\alpha$ as 
\begin{equation}
 \alpha_{2 L}(t)= \begin{cases}\frac{t+2 L}{L} \alpha(-L)+\frac{t+L}{L}\frac{\dist_N}{2} & \text { if }-2 L \leq t \leq-L, \\ \alpha(t)& \text { if }-L \leq t \leq L, \\ \frac{2 L-t}{L} \alpha(L)+\frac{t-L}{L}\frac{\dist_N}{2} & \text { if } L \leq t \leq 2 L. \end{cases}\label{1.2.15}
\end{equation} 
Then, the following properties hold:
\begin{itemize}
    \item[(1)]  The function $\alpha_{2 L}$ is monotonically increasing, and
$$
-\frac{\dist_N}{2}<\alpha_{2 L}(t)<\frac{\dist_N}{2}\quad (-2L<t<2L), \quad \alpha_{2 L}(\pm 2 L)= \pm\frac{\dist_N}{2}.
$$
\item[(2)] There exist constants $C_1,C_2>0$ such that 
    \begin{equation}
 \begin{cases} 
\left|\alpha_{2L}^{\prime}(t)\right|+\left|\alpha_{2L}(t)+\frac{\dist_N}{2}\right| \leq C_{1} e^{C_{2} t} &\text { when } t <0, \\
\left|\alpha_{2L}^{\prime}(t)\right|+\left|\alpha_{2L}(t)-\frac{\dist_N}{2}\right| \leq C_{1} e^{-C_{2} t} &\text { when } t>0.
\end{cases}
\end{equation}
\item[(3)] There exist constants $c_1,c_2>0$ such that
\begin{align}
\int_{-2 L}^{2 L}\left(\left|\alpha_{2 L}^{\prime}(t)\right|^{2}+\widetilde{F}\b(\alpha_{2 L}(t)\b)\right) d t &\leq c_F+c_{2} e^{-c_{1} L},\label{1.2.16}\\
\max_{t\in \br}\bigg||\alpha^\prime_{2L}(t)|-\sqrt{2\widetilde{F}\b(\alpha_{2L}(t)\b)}\bigg|&\leq c_2e^{-c_1L}.\label{1.2.19}
\end{align}

\label{lemma alpha properties}
\end{itemize}
\end{lemma}
\begin{proof}
     The  properties of $\alpha_{2 L}$ stated above follow directly from~\eqref{1.2.10-1}-\eqref{1.2.10}. For details, refer to \cite[Proposition A.4]{lin2012phase}.
\end{proof}

To simplify the notation,  denote $\Gamma:=\Gamma(0)$,  $\o^\pm:=\o^\pm_0$, $\dd_\Gamma:=\pm\dist(x,\Gamma),x\in\o^\pm$, and $\Gamma(0)\cap\p\o=\{{q^\pm}\}$. We extend $\Gamma(0)$ locally to $\mathbb{R}^n\setminus\o$ 
by elongating the tangents of $\Gamma$ at $q^\pm$, and denote the extended curve as $\widetilde{\Gamma}$. Then, there exists a constant $\delta\in(0,1)$ such that the signed distance function $\dg$, given by
\begin{equation}
    \dd_{\widetilde{
    \Gamma
    }}(x)=\begin{cases}
        \dd(x,\widetilde{\Gamma})&\text{ if } x\in \overline{\o^+},\\
        -\dd(x,\widetilde{\Gamma})&\text{ if } x\in\overline{\o^-},
    \end{cases}
\end{equation}
 is $C^1$ in $\widetilde{\Gamma}(\delta):=\{x\in\overline{\o}:\dg<\delta\}\cap\overline{\o}$ and satisfies (due to the  non-tangential intersections with $\p\o$)
\begin{equation}\label{equ }
    |\dg(x)|\lesssim|\dd_\Gamma(x)| \text{ for }x\in\widetilde{\Gamma}(\delta).
\end{equation}
Let $p^{ \pm} \in N^{ \pm} \text { with }\left|p^{+}-p^{-}\right|=\dist_N$. Fix a constant $\beta\in(0,1)$, and define the initial data $\uez:\o\to\rk$ as
    \begin{equation}
\uez(x)=\frac{p^{+}+p^{-}}{2}+\alpha_{\e^{-\beta}}\dge \frac{p^{+}-p^{-}}{\left|p^{+}-p^{-}\right|},\label{1.2.14}
\end{equation}
where $\alpha_{\e^{-\beta}}$ is obtained by setting $2L=\e^{-\beta}$ in~\eqref{1.2.15}.  Clearly, $$\|\uu_{\e,0}\|_{L^\infty(\o)}\leq R_N,$$
where $R_N$ is given in~\eqref{R_0}.
Let $\alpha_\e:=\alpha_{\e^{-\beta}}$ for simplicity. From the definition, we obtain
\begin{equation}\label{app cal 1}
    |\nabla\uez|=\frac{1}{\e}\alpha_\e^{\prime}\dge,\quad F(\uez)=\widetilde{F}\lef\alpha_\e\dge\rig,
\end{equation}
and
\begin{align}\label{vaule}
    \uez(x)=  p^\pm\quad\text{ for }x\in\o^\pm\setminus\widetilde{\Gamma}(\e^\beta).
\end{align}
For \(x\in\o^-\), from the definition of $\df$ (cf.~\eqref{df}), we have
\begin{align}
\dd_F(\uez)=\int^{\dd_N(\uez)}_0\sqrt{2f(\lambda)}d\lambda
    =&\int^{\frac{\dist_N}{2}+\alpha_\e\dge}_0\sqrt{2f(\lambda)}d\lambda\\
    =&\int_{-\frac{\dist_N}{2}}^{\alpha_\e\dge}\sqrt{2\widetilde{F}(\lambda)}d\lambda.\label{df es1}
     \end{align}
Similarly, for $x\in\o^+$,
\begin{equation}
   \cf- \dd_F(\uez)=\int^{\frac{\dist_N}{2}}_{\alpha_\e(\frac{\dg}{\e})}\sqrt{2\widetilde{F}(\lambda)}d\lambda.\label{df es2}
\end{equation}
\vspace{0.5em}
\noindent\textbf{Energy decomposition and estimates.}  
We decompose the relative energy $\E(0)$ as follows,
\begin{equation}
\begin{aligned}
    \E(0)=&\int_\o\left(\frac{\e}{2}\left|\nabla\uez\right|^2+\frac{1}{\e}F(\uez)-\bxi\cdot\nabla\b(\df(\uez)\b)\right)dx\\&+\int_{\p\o}\big(\sigma(\uez)-\dd_F(\uez)\cos\alpha \big) d\hn\\
    :=&\mathbf{i}_\e+\mathbf{ii}_\e,
\end{aligned} 
\end{equation}
and recall that
\begin{equation}
    \B(0)=\int_\o\b( c_F\chi_{\o^+}-\dd_F(\uez)\b)\vartheta(\cdot,0) dx.
\end{equation}

First, we estimate \(\mathbf{i}_\varepsilon\). Using~\eqref{app cal 1}-\eqref{df es2} along with the estimate~\eqref{1.2.19}, we have
\begin{equation}
  \begin{aligned}
\mathbf{i}_\e
 =&\int_\o\frac{1}{2\e}\lef\alpha_\e^{\prime}\dge\rig^2+\frac{1}{\e}\widetilde{F}\lef\alpha_\e\dge\rig-\bxi\cdot\nabla\lef\dd_F(\uez)\rig dx\\
    \leq&\int_\o\frac{1}{\e}\lef\alpha_\e^\prime\dge\rig^2-\frac{1}{\e}\bxi\cdot\nabla\dg\alpha_\e^\prime\dge\sqrt{2\widetilde{F}\lef\alpha_\e\dge\rig}dx\\
    &+\int_\o \frac{1}{\e
    } \widetilde{F}\lef\alpha_\e\dge\rig dx- \int_\o\frac{1}{2\e}\lef\alpha_\e^{\prime}\dge\rig^2dx
    \\
    \overset{\eqref{1.2.19}}\leq&\int_\o\frac{1}{\e}\lef1-\bxi\cdot\nabla\dg\rig\lef\alpha_\e^\prime\dge\rig^2dx+\frac{C}{\e}e^{-c_2\e^{\beta-1}}\\
    =&\int_{\widetilde{\Gamma}(\e^\beta)}\frac{1}{\e}\lef1-\bxi\cdot\nabla\dg\rig\lef\alpha_\e^\prime\dge\rig^2dx+\frac{C}{\e}e^{-c_2\e^{\beta-1}}.
\end{aligned}  
\end{equation}
Note that,  $\frac{C}{\e}e^{-c_2\e^{\beta-1}}$ is higher-order infinitesimal of $\e^2$ as $\e\downarrow0$. Based on the regularity of $\bxi$  given in~\eqref{eq_regularityXi} and the condition~\eqref{eq_consistencyProperty}, by applying a Taylor expansion and using \eqref{equ }, we deduce that
\begin{equation}
    |1-\bxi(x,0)\cdot\nabla\dg(x)| \lesssim\dd^2_\Gamma(x)\lesssim\dg^2(x)\quad \text{ for }x\in\widetilde{\Gamma}(\delta).
\end{equation}
Then, by applying both the Federer's co-area formula and the formula of change of variables, we obtain
\begin{align*}
    &\int_{\widetilde{\Gamma}(\e^\beta)}\frac{1}{\e}\lef1-\bxi\cdot\nabla\dg\rig\lef\alpha_\e^\prime\dge\rig^2dx\\
    \lesssim&\int_{\widetilde{\Gamma}(\e^\beta)}\frac{1}{\e}\dd_{\widetilde{\Gamma}}^2\lef\alpha_\e^\prime\dge\rig^2dx\\
     =&\int^{\e^\beta}_{-\e^\beta}\frac{1}{\e}s^2\lef\alpha_\e^\prime(\frac{s}{\e})\rig^2\times\hn\{\dg=s\}ds\\
    \lesssim&\e^2\int_{-\e^{\beta-1}}^{\e^{\beta-1}}s^2\lef\alpha_\e^\prime\lef s\rig\rig^2ds\overset{\eqref{1.2.16}}\lesssim \e^2.
\end{align*}


Next, we estimate $\mathbf{ii}_\e$. Define $g:=\sigma-\df\cos\alpha:\rk\to \br^+$. Then $g(N)=\{0\}$ and $g$ is Lipschitz continuous due to the Lipschitz continuity of $\df$ and $\sigma$.  We have
\begin{align}
     |g(\uez)|\leq\begin{cases}
    Lip(g) |\alpha_\e\dge+\frac{\dist_N}{2}|\quad\text{ if }x\in\p\o\cap\overline{\widetilde{\Gamma}(\delta)}\cap\{\dg<0\},\\
     Lip(g)|\alpha_\e\dge-\frac{\dist_N}{2}|\quad\text{ if }x\in\p\o\cap\overline{\widetilde{\Gamma}(\delta)}\cap\{\dg>0\}.
    \end{cases}
\end{align}
Using~\eqref{vaule}, and applying the formula of change of
variables along with a scaling argument, we get
\begin{align*}
     \int_{\p\o} |g(\uez)|\hn
     \lesssim&\int_{-\e^\beta}^{0}\bigg|\alpha_\e(\frac{s}{\e})+\frac{\dist_N}{2}\bigg|ds+\int_0^{\e^\beta}\bigg|\alpha_\e(\frac{s}{\e})-\frac{\dist_N}{2}\bigg|ds\\
     \lesssim&\e\int_{-\e^{\beta-1}}^0\bigg|\alpha_\e(s)+\frac{\dist_N}{2}\bigg|ds+\e\int_0^{\e^{\beta-1}}\bigg|\alpha_\e(s)-\frac{\dist_N}{2}\bigg|ds\\
     \overset{\eqref{1.2.16}}\lesssim&\e.
\end{align*}

Finally, we estimate $\B(0)$. Using~\eqref{vaule},~\eqref{df es1} and~\eqref{df es2}, it follows
\begin{align}
    B_\e[\uez|\Gamma]=&\left\{\int_{\o^+\cap\widetilde{\Gamma}(\e^\beta)}+\int_{\o^-\cap\widetilde{\Gamma}(\e^\beta)}\right\}\b(c_F\chi_{\o^+}-\dd_F(\uez)\b)\vartheta (\cdot,0)dx\\
    \lesssim&\int_{\o^+\cap\widetilde{\Gamma}(\e^\beta)}\bigg(\left|\vartheta (\cdot,0)\right|\int^{\frac{\dist_N}{2}}_{\alpha_\e\dge}\sqrt{2\widetilde{F}(\lambda)}\bigg) dx\\
    &+\int_{\o^-\cap\widetilde{\Gamma}(\e^\beta)}\bigg(\left|\vartheta (\cdot,0)\right|\int_{-\frac{\dist_N}{2}}^{\alpha_\e\dge}\sqrt{2\widetilde{F}(\lambda)}\bigg) dx. 
\end{align}
From the condition~\eqref{2.16d} and~\eqref{equ }, we have
\begin{equation}
    \left|\vartheta(x,0)\right| \lesssim\dd_\Gamma(x)\lesssim\dg(x)\quad\text{ for }x\in\widetilde{\Gamma}(\delta).
\end{equation}
By both the Federer’s co-area formula and the formula of change of
variables, using the  local boundedness of $\widetilde{F}$, we get
\begin{equation}
    \begin{aligned}  B_\e[\uez|\Gamma]\lesssim&\int_{\o^+\cap\widetilde{\Gamma}(\e^\beta)}\left|\dg\right|\left|\frac{\dd_N}{2}-\alpha_\e\dge\right|dx\\
    &+\int_{\o^-\cap\widetilde{\Gamma}(\e^\beta)}\left|\dg\right|\left|\frac{\dd_N}{2}+\alpha_\e\dge\right|dx\\
    \lesssim&\e^2\int_{-\e^{\beta-1}}^0s\bigg|\alpha_\e(s)+\frac{\dist_N}{2}\bigg|ds+\e^2\int_0^{\e^{\beta-1}}s\bigg|\alpha_\e(s)-\frac{\dist_N}{2}\bigg|ds\\
    \overset{\eqref{1.2.16}}\lesssim& \e^2.
    \end{aligned}
\end{equation}

\end{document}